\documentclass[a4paper]{article}
\usepackage{amssymb}
\usepackage{latexsym}
\newtheorem{theorem}{Theorem}
\newtheorem{definition}{Definition}
\newtheorem{propo}{Proposition}

\def\Proof{\noindent{\it Proof}\rm}
\def\qed{\vbox
{\hrule\hbox{\vrule\hbox to 5pt{\vbox to
5pt{\vfil}\hfil}\vrule}\hrule}}
\def\endproof{\unskip \nobreak \hskip0pt plus 1fill
\qquad \qed \par}

\font\tenscr=rsfs10 

\font\sevenscr=rsfs7 

\font\fivescr=rsfs5 
\skewchar\tenscr='177 \skewchar\sevenscr='177
\skewchar\fivescr='177
\newfam\scrfam \textfont\scrfam=\tenscr \scriptfont\scrfam=\sevenscr
\scriptscriptfont\scrfam=\fivescr
\def\scr{\fam\scrfam}
\def\ca{{\scr A}}
\def\cj{{\scr J}}
\def\ct{{\scr T}}

\font\tencmmib=cmmib10 \skewchar\tencmmib='177
\let\bsf=\tencmmib

\def\fg{\hbox{\bsf g}}

\def\fs{\hbox{\bsf s}}
\def\fe{\hbox{\bsf e}}
\def\ft{\hbox{\bsf t}}
\def\fp{\hbox{\bsf p}}

\def\Cl{{\cal C}\mskip-4mu \ell}

\def\C{{\mathbb{C}}}
\def\H{{\mathbb{H}}}
\def\K{{\mathbb{K}}}
\def\N{{\mathbb{N}}}
\def\R{{\mathbb{R}}}
\def\Z{{\mathbb{Z}}}

\let\a=\ast
\let\e=\varepsilon
\let\vi=\varphi
\let\p=\prime
\let\l=\langle
\let\r=\rangle
\let\w=\wedge
\let\d=\partial
\let\b=\bar
\def\Ad{\mathop{\rm Ad}\nolimits}
\def\Aut{\mathop{\rm Aut}\nolimits}
\def\coker{\mathop{\rm Coker\,}\nolimits}
\def\curv{\hbox{\rm curv }}
\def\End{\mathop{\rm End\,}\nolimits}
\def\HOM{\mathop{\rm HOM\,}\nolimits}
\def\Hom{\mathop{\rm Hom\,}\nolimits}
\def\id{\mathop{\rm id\,}\nolimits}
\def\im{\mathop{\rm Im\,}\nolimits}
\def\ind{\mathop{\rm ind\,}\nolimits}
\def\Iso{\mathop{\rm Iso\,}\nolimits}
\def\Ker{\mathop{\rm ker\,}\nolimits}
\def\mod{\mathop{\rm mod\,}\nolimits}
\def\Span{\mathop{\rm span\,}\nolimits}
\def\sgn{\mathop{\rm sgn\,}\nolimits}
\def\tr{\mathop{\rm tr\,}\nolimits}

\def\inp{\,\hbox{\vrule height5pt width0.35pt
    {\vrule height0.35pt width6.5pt}}\,}
\def\L{\mathop{\textstyle\bigwedge}\nolimits}

\def\oco{\mathop{\otimes\cdots\otimes}}
\def\moco{\mathop{\;\otimes\cdots\otimes\;}}

\def\o#1{\mathop{\overline{#1}}}
\def\ctr#1{\hfil#1\hfil}
\def\vl{\vrule height 17truept depth 7truept}
\def\cancel#1#2{\ooalign{$\hfil#1\mkern1mu/\hfil$\crcr$#1#2$}}
\def\nD{\mathrel{\mathpalette\cancel\partial}}
\def\nd{\mathrel{\mathpalette\cancel D}}
\def\vl{\vrule height 17truept depth 7truept}

\begin{document}
\centerline{\Large\bf {On the definition of geometric Dirac
operators}}
\medskip
\medskip
\medskip
{\centerline{\it Herbert Schr\"{o}der}}
\bigskip
\bigskip
\bigskip
\bigskip
\bigskip
\bigskip
\bigskip

For the definition of a spin$^c$ structure and its associated
Dirac operators there can be found two different approaches in the
literature. One of them uses lifts of the orthonormal frame bundle
to principal spin$^c$ bundles (cf.\ [Gil], [GH], [Frie] or [LM]) and
the other one irreducible
representations of the complex Clifford bundle (cf.\ [BD] or [Kar1,2]).
The first approach
is an offspring of vector and tensor calculus in its modern form
as shaped by E.\ Cartan and Ch.\ Ehresmann whereas the second
approach is rooted in physics, in particular, relativistic quantum
mechanics. Although the second approach is favored nowadays, in
defining spin structures most authors still rely on the first
method. In this expository note we give a definition of spin
structure and the corresponding Spin-Dirac operator purely in the
spirit of irreducible representations and prove its equivalence
with the usual definition.

This seems to be well known to people working in noncommutative
geometry. At least it is used and taken for granted e.g.\ in
[Con1-3], [Ren], and [Var]. The purpose of our note is to make
this method accessible to a wider audience in mathematics and in
physics and to direct attention to the so far mostly ignored work
of G.\ Karrer who introduced spin$^c$-structures in this way
already in 1962 and published his results in 1963 [Kar1] and 1973
[Kar2] (unfortunately in German); usually, this approach is
credited to A.\ Connes (cf.\ [BD]).

Spinors first appeared in the theory of representations of the
orthogonal group, in fact of its Lie algebra, in 1913 [Car] and
then again in 1927 in connection with the Dirac equation [Dir].
The Schr\"{o}dinger equation of classical quantum mechanics
\[\frac{1}{i}\frac{\partial}{\partial t}\psi+\Delta\psi=0 \]
(without external electro-magnetic field and ripped off of any
physical meaning by setting the usual constants $\hbar$, $m$ and
$c$ equal to 1) is of first order in the time variable and
invariant under Galilei transformations, i.e.\ time and spatial
translations and spatial rotations. To get a relativistic analogue
again of first order in $t$ and invariant under Lorentz
transformations, P.A.M.\ Dirac was looking for a square root of
the d'Alembert operator $\Box = \frac{\partial^2}{\partial
x_0}-\sum_{j=1}^3 \frac{\partial^2}{\partial x_j^2}$ which governs
the Klein-Gordon equation. He found
\[\Box~ = ~\nD^2 ~= ~\Big(\sum_{j=0}^3 A_j {\d\over\d x_j}\Big)^2,\ \]
where
\[A_0 = \left(\begin{array}{cc} I_2 & 0 \\
0 & -I_2\\
\end{array}\right)\quad \mbox{and} \quad A_j =
\left(\begin{array}{cc} 0 & \sigma_j\\ -\sigma_j &
0\\\end{array}\right)\,
\in
M_2\big(M_2(\C)\big),\quad j=1,2,3\,,\]
with the Pauli matrices
\[\sigma_1 = \left(\begin{array}{cc}
0&1\\ 1&0\\
\end{array}\right) \,,\quad \sigma_2 =
\left(\begin{array}{cc}0&-i\\ i&0\\ \end{array}\right)\quad \hbox{
and } \quad \sigma_3=\left(\begin{array}{cc}1&0\\
0&-1\\\end{array}\right).\]
The physicists encountered however
difficulties in combining Dirac's first order equation $\nD\psi=0$
for the relativistic electron with the needs of general
relativity, since the spinors $\psi$ did not transform like
vectors or tensors, and so, in first instance, had no geometrical
meaning. We quote the words of C.G.\ Darwin in 1928 [Dar]:
\begin{quote}
\baselineskip10pt {\footnotesize The relativity theory is based on
nothing but the idea of invariance and develops from it the
conception of tensors as a matter of necessity; and it is rather
disconcerting to find that apparently something has slipped
through the net, so that physical quantities exist, which it would
be, to say the least, very artificial and inconvenient to express
as tensors.}
\end{quote}
The following years saw various attempts to find a non-local version of the 
Dirac operator, which had to act on spinors; cf.\ [vdW1] for a detailed 
historical survey of concept of spin in physics. But still in 1937 E.\ Cartan 
in his book ``La Th\'eorie des Spineurs'' noted unsurmountable difficulties 
to apply techniques of classical tensor calculus to spinors. Only in the 
fifties with the invention of principal bundles and its connections the 
spinors found their appropriate place in Riemannian geometry. It became 
possible to define the covariant derivative of spinors and finally around 
1960 to define the Dirac operator. This has been achieved by E.\ Kaehler for 
the Dirac operator $d+\delta$ in 1961 [Kae] and by M.F.\ Atiyah and I.M.\ 
Singer for the Spin-Dirac operator in 1962 [AS]. 

The definition of spin structures consists of two parts, a local
and a global one. The local part is purely algebraic and will be
treated in the first two sections. Here we sketch the most
important results concerning Clifford algebras and refer to e.g.\
[Che], [Krb] or [LM] for a more detailed account. The global part
is of topological nature. It will be exposed in the third section.
In the fourth section we discuss spin structures in the setting of
principle bundles. We conclude with the definition and some
elementary properties of some geometric Dirac operators.
\bigskip
\bigskip
\section{\bf Clifford algebras}
\bigskip

Clifford algebras solve an algebraic existence problem. To see
this recall that the field of complex numbers arises in two ways.
In the first instance it is merely a vector space that helps
parametrize the Euclidean plane $\R^2$ but in the second it is an
algebra extending the real number field in which square roots
exist and which contains an image of the group of rotations. In
particular, only by this property we comprehend the law of
multiplication of two negative numbers: $(-1)(-1)=1$, since
$-1=i^2$ is the composition of two rotations by $90$ degrees. As
is well known it took R.W.\ Hamilton ten years to find out in 1843
that there is no analogue in 3-space. One has to step out of
ordinary space to find an algebra which contains $\R^3$ as well
its rotations, viz.\ the skew field of quaternions. What is the
appropriate generalization to arbitrary dimension? Starting from a
real vector space $E$, one has e.g.\ the exterior algebra $\L E$
introduced by H.G.\ Grassmann in 1844. It contains $E$ and its
multiplication $\w$ is anti-commutative on basis vectors:
\[e_i\w e_j + e_j\w e_i = 0\,.\]
But then basis vectors $e_i$ are nilpotent, $e_i\w e_i = 0$. What we really 
need is a new multiplication $\cdot$ such that the basic vectors satisfy $e_i 
\cdot e_i = -1$. How to come to terms with this has first been observed by 
W.K.\ Clifford in 1876: 
\begin{quote}
\baselineskip10pt {\footnotesize The system of quaternions differs
from this, first in that the squares of the units, instead of
being zero, are made equal to $-1$; and secondly in that the
ternary product $\iota_1\iota_2\iota_3$ is made equal to
$-1$.\dots\\ I shall now examine the consequence of making, in a
system of $n$ alternate numbers $\iota_1,\iota_2,\dots,\iota_n$,
the {\it first} of the modifications just named; namely I shall
suppose that the square of each of the units is $-1$.}
\end{quote}
After the advent of modern abstract algebra the construction of Clifford's 
``geometric algebra'' runs as follows. We choose an inner product 
$\l\cdot,\cdot\r $ on $E$, i.e., we assume a Euclidean vector space 
$(E,\l\cdot,\cdot\r )$, and with respect to this inner product we choose an 
orthonormal basis $(e_i)_{1\le i\le n}$ that satisfies 
\[e_i \cdot e_j + e_j \cdot e_i = - 2\delta_{ij} =-2 \l e_i,e_j\r
\,.\]
This obtains from assuming
\[v\cdot v + \l v,v\r \,= 0\,\]
for any $v\in E$. Just like the exterior algebra the new algebra
we are looking for can now be constructed as a quotient of the
tensor algebra $\ct(E)$. Here we have to consider the two-sides
ideal $\cj(E)\subset\ct(E)$, which is generated by elements
$v\otimes v + \l v,v\r 1$, $v\in E$.

\begin{definition}\quad The $\R$-algebra $\Cl(E) =
\ct(E)/\cj(E)$ (corresponding to a given Euclidean structure) is
called the Clifford algebra of $E$. In case of $E=\R^n$ with its
standard Euclidean structure we write $\Cl_n = \Cl(\R^n)$.
\end{definition}

The product in $\Cl(E)$ will be denoted by $\cdot$, i.e.\ for $u,v\in \Cl(E)$ 
with $u=\pi(\tilde u)$, $v=\pi({\tilde v})$, where $\pi : \ct(E) \to \Cl(E)$ 
denotes the natural projection, let $u\cdot v=\tilde u \otimes {\tilde v} + 
\cj(E)$. We also denote by $\iota_E:E\to\Cl(E)$ the restriction of $\pi$ to 
$E$. 

Just like the tensor algebra and the exterior algebra the Clifford
algebra solves a universal problem.

\begin{theorem} Given an associative unital $\R$-algebra $A$ (with
unit 1) and a linear map $f:E\to A$ with $f(v) \cdot f(v) = -\l v,v\r 1$ for 
all $v\in E$, there is a unique homomorphism of $\R$-algebras, $\tilde 
f:\Cl(E) \to A$,  such that the following diagram commutes 

\unitlength10mm 
\begin{picture}(6,4)(-4.5,-1)
\put(0,0){E} \put(0.5,0.1){\vector(1,0){3.3}} \put(4,0){A}
\put(0.1,0.5){\vector(0,1){1.5}} \put(-0.2,2.3){$\Cl(E)$}
\put(0.5,2){\vector(2,-1){3.3}}
\put(-0.5,1.3){$\scriptstyle{\iota_E}$}
\put(1.7,0.4){$\scriptstyle{f}$}
\put(1.7,1.7){$\scriptstyle{\tilde f}$}
\end{picture}
\\
In particular, the algebra $\Cl(E)$ together with the map $\iota_E : E\to 
\Cl(E)$ satisfying $\iota_E(v)^2 = -\l v,v\r 1$ is uniquely determined by 
this property up to isomorphism. 
\end{theorem}

\Proof: We have
\[\iota_E(v)^2 = \pi(v)^2 = v\otimes v + \cj(E) =
-\l v,v\r 1 + \cj(E) = -\l v,v\r 1\quad {\hbox{ for }} v\in E\]
(here $1=1+\cj(E)$ is the unity of $\Cl(E)$). Since $\ct(E)$ is
generated by $E$ as an algebra and since $\pi$ is surjective,
$\Cl(E)$ is generated by $\iota_E(E) $. Now given a linear map
$f:E\to A$ with $f(v)^2 = -\l v,v\r 1_A$, $v\in E$, we have an
extension to a homomorphism of algebras, $\otimes f : \ct(E) \to
A$, given by
\[\otimes f\big(v\otimes v+\l v,v\r 1\big) = f(v)^2
+\l v,v\r 1_A = 0\,,\] and hence factorizes to a homomorphism of
algebras, ${\tilde f} : \Cl(E) \to A$. Then for $v\in E$ we have
\[{\tilde f} \circ \iota_E(v) = {\tilde f}\circ \pi(v) = \otimes
f(v) = f(v)\]
and ${\tilde f}$ is uniquely determined since
$\iota_E(E)$ generates $\Cl(E)$.
\endproof
\medskip

Clifford algebras have entered quite different branches of modern
mathematics and physics in the 100 years since their introduction
by W.K.\ Clifford in 1876 [Cli] and independently by R.\ Lipschitz
in 1880 [Lip]; cf.\ also his letter from Hades written by his
medium A.\ Weil [Wei]. Clifford's main purpose was to generalize
H.G.\ Grassmann's exterior algebra and R.W.\ Hamilton's
quaternions, whereas Lipschitz was looking for a parametrization
of orthogonal transformations of $\R^n$. That Clifford algebras
indeed meet both purposes turned out in 1935, when R.\ Brauer and
H.\ Weyl [BW] gave a very elegant representation of the spin
group.

In 1954 C.\ Chevalley [Che] gave the concise construction
presented above. It allows the inner product to be replaced by an
arbitrary symmetric bilinear form $\sigma : E\times E\to \K$, or,
more precisely, by the corresponding quadratic form $Q$, and
$\K=\R$ or $\C$ by any field. We preferably consider $\K=\R$ or
$\C$ depending on $E$ being a real or a complex vector space. In
general, one obtains Clifford algebras $\Cl(E,Q)$, in particular,
for $Q=0$ the exterior algebra. On $E=\R^{r+s}$ one considers the
quadratic forms
\[Q_{r,s}(x) = \sum_{i=1}^r x_i^2 - \sum_{i=r+1}^{r+s} x_i^2\]
yielding the Clifford algebras $\Cl_{r,s}$. We take a look at some special 
examples. 
\medskip

\noindent
{\bf{Examples}}\quad 1. For $E=\R$ with inner product $\l x,y\r  =
x y$ we have $\Cl(\R) = \Cl_1 = \C$. For if $\iota_{\R}(x) = ix$,
$x\in \R$, the algebra $\C$ is generated by $\iota_{\R}(\R)$ since
$\iota_{\R}(x)^2 = -\l x,x\r 1$. Given an algebra $A$ and $f:\R
\to A$ as above with $f(x)^2 = -\l x,x\r 1_A$, we get
\[f(x) = x f(1)\,,\]
since $f$ is linear, and
\[f(1)^2=-1_A\,.\]
Defining
\[{\tilde f}(x+iy)=x1_A + y f(1),\quad x,y\in \R,\]
we obtain a homomorphism of $\R$-algebras and
\[{\tilde f}\circ
\iota_{\R} (y) = {\tilde f}(i y) = y f(1) = f(y),\quad y\in \R.\] \noindent 
2. The Clifford algebra $\Cl_2$ is isomorphic with the skew field of 
quaternions, $\H$, which is generated by $i\iota_{\R^2}(e_1)$ and 
$j=\iota_{\R^2}(e_2)$, since $k=ij$ and $i^2=j^2=k^2=-1$, if $\{e_1,e_2\}$ 
denotes the standard basis of $\R^2$. 
\medskip

\noindent {\bf{Remarks}} 1. By the universal property any isometry
$f:\big(E,\l \cdot,\cdot\r\big) \to \big(E^\p,\l\cdot,\cdot\r ^\p\big)$
induces a homomorphism of algebras $\Cl(f) : \Cl(E) \to \Cl(E^\p)$: One
simply has to lift the map $\bar f=\iota_{E'}\circ f$ that satisfies

\[\bar f(v)^2 =\iota_{E^\p}\big(f(v)\big)^2 =
-\,\l f(v),f(v)\r ^\p 1=-\,\l v,v\r 1 \]

as in the (commutative) diagram

\unitlength1cm
\begin{picture}(6,4)(-4.5,-1)
\put(0,0){$E$}
\put(0.5,0.1){\vector(1,0){3.3}}
\put(4,0){$E'$}
\put(0.1,0.5){\vector(0,1){1.5}}
\put(-0.3,2.3){$\Cl(E)$}
\put(4.2,0.5){\vector(0,1){1.5}}
\put(3.9,2.3){$\Cl(E')$}
\put(0.5,0.3){\vector(2,1){3.3}}
\put(0.6,2.4){\vector(1,0){3.1}}
\put(-0.5,1.3){$\scriptstyle{\iota_E}$}
\put(1.9,0.4){$\scriptstyle{f}$}
\put(1.9,1.5){$\scriptstyle{\tilde f}$}
\put(4.5,1.3){$\scriptstyle{\iota_{E'}}$}
\put(1.7,2.6){$\scriptstyle{\Cl(f)}$}
\end{picture}

\noindent
Given another isometry $g :\big(E^\p,\l\cdot,\cdot\r ^\p\big) \to
\big(E^{\p\p},\l\cdot,\cdot\r ^{\p\p}\big)$, one has
\[\Cl(g\circ f) = \Cl(g) \circ \Cl(f).\]
Therefore, $\Cl:O(E) \to \Aut~\Cl(E)$ defines a homomorphism of
groups.\\
2. The involution $\alpha:E\to E$, $\alpha(v) = -v$,
$v\in E$, extends to an involution of $\Cl(E)$ again denoted by
$\alpha$. Using $\alpha$ one defines a $\Z_2$-grading
\[\Cl(E) = \Cl(E)^0 \oplus \Cl(E)^1\]
by $\alpha\big|_{\Cl(E)^j} = (-1)^j \id$, $j=0,1$; since $\alpha$
is a homomorphism, we have
\[\Cl(E)^i \cdot \Cl(E)^j \subset \Cl(E)^{i+j~ {\rm mod}~2}\]
turning $\Cl(E)^0$ into a subalgebra.

\begin{propo} Given $v,w\in E$ with $\l v,w\r =0$
one has
\[\iota_E(v) \cdot \iota_E(w) + \iota_E(w) \cdot \iota_E(v) = 0\,.\]
More generally, given $x=\prod_{\ell=1}^n \iota_E(v_\ell)\in
\Cl(E)^i$ and $y=\prod_{k=1}^m \iota_E(w_k)\in \Cl(E)^j$ with $\l
v_\ell,w_k\r = 0$ for all $\ell$ and $k$, one has
\[x\cdot y =(-1)^{ij} y\cdot x\,.\]
\end{propo}

\Proof: We compute $\iota_E(v+w)^2$ in two ways. One the one hand
\[\iota_E(v+w)^2 = -\l v+w,v+w\r 1 = -\l v,v\r 1 - \l w,w\r 1\]
and on the other hand
\[\iota_E(v+w)^2 =  \iota_E(v)^2 +
\iota_E(w)^2 + \iota_E(v)\cdot \iota_E(w) + \iota_E(w)\cdot
\iota_E(v)\,.\]
Equating both sides gives the first assertion. The
second one follows by induction since $\iota_E(E)\subset
\Cl(E)^1$.
\endproof
\medskip

In order to prove the basic structure theorem for Clifford
algebras we need the notion of graded tensor product of two graded
algebras. Given two unital $\R$-algebras $A$ and $B$ with units
$1_A$ and $1_B$, resp., the tensor product $A\otimes B$ turns into
an $\R$-algebra if we put
\[(a\otimes b)(a^\p\otimes b^\p) = a a^\p\otimes bb^\p\,\]
for $a,a^\p\in A$, $b,b^\p\in B$. If $A$ and $B$ are
$\Z_2$-graded, i.e.\ $A=A^0\oplus A^1$ and $B=B^0\oplus B^1$ a
$\Z_2$-grading of $A\otimes B$ is defined by
\begin{eqnarray*}
(A\otimes B)^0 & = & A^0\otimes B^0\oplus A^1\otimes B^1\\ 
(A\otimes B)^1 & = 
& A^1\otimes B^0\oplus A^0\otimes B^1, 
\end{eqnarray*}

\noindent
where the product is now given by
\[(a\otimes b)(a^\p \otimes b^\p) = (-1)^{ij} aa^\p\otimes bb^\p\]
for $a^\p\in A^i$, $b\in B^j$. To distinguish the two tensor products, we 
denote the graded tensor product of $A$ and $B$ by $A\hat{\otimes} B$. 

\begin{theorem} Any orthogonal splitting $E=E_1\oplus E_2$ gives rise to a 
canonical isomorphism of algebras $\Cl(E)$ and $\Cl(E_1) \hat{\otimes} 
\Cl(E_2)$. 
\end{theorem}

\Proof: We start with $f:E\to\Cl(E_1)\hat{\otimes}\Cl(E_2)$ defined by 
\[f(v_1+v_2) = \iota_{E_1}(v_1)\otimes 1 + 1\otimes
\iota_{E_2}(v_2)\,,v_k\in E_k.\] Since $\iota_{E_k}(v_k) \in
\Cl(E_k)^1$, $1\in \Cl(E_k)^0$, and $v_1\bot v_2$, we get
\begin{eqnarray*}
f(v_1+v_2)^2 &=& \iota_{E_1}(v_1)^2 \otimes 1 + 1\otimes \iota_{E_2}(v_2)^2\\
&=& \big(-\l v_1,v_1\r - \l v_2,v_2\r \big) 1\otimes 1\\
&=& -\l
v_1+v_2,v_1+v_2\r 1\otimes 1\,
\end{eqnarray*}

\noindent
hence a unique homomorphism ${\tilde f}:\Cl(E) \to \Cl(E_1)
\hat\otimes \Cl(E_2)\,$ by Theorem 1. Likewise the isometries $i_1
: E_1\to E$ and $i_2:E_2\to E$ induce homomorphisms $\Cl(i_k)$,
$k=1,2$, and for $x\in \Cl(E_1)^i$, $y\in \Cl(E_2)^j$ one has
\[\Cl(i_1)(x)\cdot \Cl(i_2)(y) = (-1)^{ij} \Cl(i_2)(y)\cdot
\Cl(i_1)(x)\]
by the Proposition. Hence $\tilde g : \Cl(E_1)
\hat\otimes \Cl(E_2) \to \Cl(E)$ defined by
\[\tilde g(x\otimes y)
= \Cl(i_1)(x) \cdot \Cl(i_2)(y)\,,\quad x\in \Cl(E_1),\, y\in
\Cl(E_2),\]
is a homomorphism; and a straightforward computation
on generators shows that $\tilde f$ and $\tilde g$ are mutual
inverses.
\endproof

{\bf Remark}\quad An analogous result holds in case of a direct composition
$E=E_1\oplus E_2$ into $\K$-vector spaces with respect to a quadratic
form $Q=Q_1\oplus Q_2$.
\medskip

\noindent
{\bf Corollary}~ {\it Given an orthonormal basis $(e_i)_{1\le i\le
n}$ of $\big(E,\l\cdot,\cdot\r\big)$ one obtains a basis
\[\{\iota_E(e_{k_1}) \cdots \iota_E(e_{k_r}) \mid 1\le k_1 <\cdots
< k_r \le n,\; r\ge 0\}\] of $\Cl(E)$. In particular, $\dim \Cl(E)
= 2^n$ and multiplication in $\Cl(E)$ is determined by the
relations
\[\iota_E(e_k)\cdot \iota_E(e_k) = -1\,,\quad
\iota_E(e_k)\cdot \iota_E(e_\ell) + \iota_E(e_\ell)\cdot
\iota_E(e_k) = 0\quad {\hbox{ for }} k\ne \ell\,.\] Moreover one
has $\Cl(E)^i=\Span\{\iota_E(e_{k_1})\cdots \iota_E(e_{k_r})\mid
r=i\mod 2\}$.}
\smallskip

\Proof: We decompose $E$ orthogonally into
\[E=\bigoplus_{k=1}^n\R e_k\]
and apply Theorem 2 repeatedly using Example 1:
\[\Cl(E)
\cong \big(\R\oplus \R \iota_E(e_1)\big) \hat\otimes \cdots
\hat\otimes \big(\R\oplus \R\, \iota_E(e_n)\big).\] It is clear
that the multiplication is determined by the given relations. From
\[\alpha\big(\iota_E(e_{k_1})\cdots
\iota_E(e_{k_r})\big)=(-1)^r\iota_E(e_{k_1})\cdots
\iota_E(e_{k_r})\]
we obtain the final assertion.
\endproof
\medskip

From the Corollary we see that $\iota_E : E\to \Cl(E)$ is
injective. Therefore, we can identify $E$ with its image
$\iota_E(E)$ and multiply $v,w\in E$ within $\Cl(E)$, i.e., we
write $v\cdot w$ instead of $\iota_E(v)\cdot \iota_E(w)$. We also
extend the inner product of $E$ to the inner product of $\Cl(E)$
that renders the basis of the Corollary an orthonormal basis. Also
note that an isomorphism $\Cl_{n-1}\cong\Cl^0_n$ is induces by
$e_k\mapsto e_k\cdot e_n$, $k=1,\dots,n-1$, given an orthonormal
basis $\{e_1,\dots,e_n\}$ of $\R^n$.

Since $\L E$ and $\Cl(E)$ have the same dimensions they are
isomorphic as $\R$-vector spaces although not as $\R$-algebras. A
canonical homomorphism $\phi:\L E\to\Cl(E)$ is given by
\[\phi(v_1
\w \cdots \w v_k) = \frac{1}{k!} \sum_{\sigma\in S_k} (\sgn
\sigma)~ v_{\sigma(1)} \cdots v_{\sigma(k)}.\]

It is one-to-one since
\[\phi(e_{j_1} \w \cdots \w e_{j_k}) = e_{j_1} \cdots e_{j_k}\]
and actually an isometry if $\L E$ is
equipped with the appropriate inner product. The inverse
isomorphism $\sigma : \Cl(E) \to \L E$ is given by
\[\sigma(x) = c(x)1\in \L E,\quad c\in\Cl(E),\]
where $1\in\R=\L^0 E$ and where
$c:\Cl(E)\to\End(\L E)$ denotes the unique extension of the linear
map $c:E\to \End(\L E)$ defined by
\[c(v)\omega=v\wedge\omega-v\inp\omega,\quad \omega\in\L E,~ v\in
E.\] 
We already mentioned the generalized Clifford algebras $\Cl_{r,s}$. It 
is easily shown that they are generated by multiplying the standard basis 
elements $e_1, \cdots, e_{r+s}$ of $\R^{r+s}$ while respecting 
\[e_k \cdot e_\ell + e_\ell \cdot e_k
= \cases{-2 \,,&$k=\ell\le r$\cr 2\,,&$k=\ell > r$\cr
0\,,&else.\cr}\hfill(\ast)\] If $n=r+s$ is even we put $\e =
e_1\cdot \dots \cdot e_n$. Using $(\ast)$ we get
\[\e^2 =(-1)^{(n-1)+(n-2)+\cdots + 1} e_1^2 e_2^2 \cdots e_n^2
= (-1)^{\frac{n(n-1)}{2}}(-1)^r 1\] and we call $\Cl_{r,s}$ positive or 
negative if $\e^2 = +1$ or $-1$, respectively. Since the index of a quadratic 
form does not depend on the chosen basis we can speak of a positive or 
negative Clifford algebra $\Cl(E,Q)$ in case of any vector space of even 
dimension and any non-degenerate quadratic form. 

\begin{theorem} Let $E=E_1 
\oplus E_2$ and $Q = Q_1 \oplus Q_2$ with dim $E_1$ even. Then 
\[\Cl (E,Q) \cong\Cl(E_1,Q_1) \otimes \Cl(E_2,\pm Q_2)\,,\]
the sign depending
on $\Cl(E_1,Q_1)$ being positive or negative, respectively.
\end{theorem}

\Proof: For $\e = e_1 \cdot \dots \cdot e_n\in\Cl(E_1,Q_1)$, $n=\dim E_1$, 
one has 
\[\e e_i = (-1)^{n-1} e_i\e = -e_i \e,\]
i.e., $\e v = -v\e $ for any $v\in E_1 \subset \Cl(E_1)$. We
define
\[\vi : E = E_1 \oplus E_2 \to \Cl(E_1,Q_1) \otimes \Cl(E_2,\pm Q_2)\]
by
\[\vi(v_1,v_2) = v_1 \otimes 1 + \e \otimes v_2\,,\quad v_i\in E_i,\; \]
and obtain
\begin{eqnarray*}
\vi(v_1,v_2)^2 &=& v_1^2 \otimes 1 + \e^2 \otimes v_2^2 + v_1 \e \otimes v_2
+ \e v_1 \otimes v_2\\ &=& v_1^2\otimes 1\pm 1\otimes v_2^2\\
&=&-\big(Q_1(v_1) + Q_2(v_2)\big)1\otimes 1\,
\end{eqnarray*}
 since
$v_2^2=-\big(\pm Q_2(v_2)\big)=\mp Q_2(v_2)$. Using Theorem 1
(more precisely the corresponding result for an arbitrary
quadratic form) we obtain a homomorphism
\[\tilde\vi : \Cl (E,Q)
\to \Cl(E_1,Q_1) \otimes \Cl(E_2,\pm Q_2)\,.\] Since dimensions
match we are reduced to verify that $\tilde\vi$ is surjective. To
this end it suffices to show that $v_1 \otimes 1$ and $1\otimes
v_2$ belong to the image of $\tilde\vi$. But now we have
\[v_1
\otimes 1 = \tilde\vi \big(\iota_E(v_1)\big) \quad \hbox{ and }
\quad 1 \otimes v_2 =
\pm\tilde\vi\big(\iota_E(v_2)\cdot\e\big)\,,\]
which concludes the
proof.
\endproof
\medskip

\begin{propo}\quad If $\dim E$ is even and $\Cl(E,Q)$
positive, then
\[\Cl(E,Q) \cong \Cl(E,-Q)\,.\]
\end{propo}

\Proof: Employing the canonical maps $\iota_\pm:E\to\Cl(E,\pm Q)$ we put 
\[\e_\pm =\iota_\pm (e_1)\cdots \iota_\pm(e_n).\]
Now $f: E \ni v \mapsto \e_+\cdot \iota_+(v) \in \Cl(E,Q)$
satisfies
\[f(v)^2 = -\e_+^2\cdot \iota_+(v)^2 = -\big(-Q(v)\big),\,\]
hence induces a homomorphism $\tilde f : \Cl(E,-Q) \to \Cl(E,Q)$.
Moreover
\[\tilde f(\e_-\cdot\iota_-(v)) = \big(f(e_1) \cdots
f(e_n)\big)\cdot f(v) =
(-1)^{n(n-1)/2}\e_+^{n+2}\cdot\iota_+(v)=\pm \iota_+(v)\,\]
whereby $\tilde f$ is surjective, hence bijective.
\endproof
\medskip

We have already determined $\Cl_{1,0}=\C$ and $\Cl_{2,0}=\H$. It
is not difficult to see that
\begin{eqnarray*}
\Cl_{0,1} &\cong& \R \oplus \R=\R (1,1)+\R (1,-1),\\ \Cl_{0,2} &\cong&
M_2(\R) \quad\hbox{with}\quad e_1 = \pmatrix{1 &0\cr 0&-1}~\hbox{ and } ~ e_2
= \pmatrix{0&1\cr 1&0}\\ \Cl_{1,1} &\cong& M_2 (\R) \quad\hbox{with}\quad e_1
= \pmatrix{0&1\cr -1&0\cr}~\hbox{ and }~ e_2 = \pmatrix{0&1\cr 1&0\cr}\,.
\end{eqnarray*}
Combining Theorem 3 with the last Proposition we obtain the following 
complete classification of Clifford algebras $\Cl_{r,s}$. 

\begin{theorem} The Clifford algebras $\Cl_{r+n,s+n}$ and
$M_{2^n}(\Cl_{r,s})$ are isomorphic, in particular
\[\Cl_{n,n}
\cong M_{2^n}(\R)\,,\quad \Cl_{r,s} \cong \cases{M_{2^s}(\Cl_{r
-s,0})\,,&$r>s$\cr M_{2^r}(\Cl_{0,s-r})\,,&$r<s$\cr}\]
and
\[\Cl_{r+8,s} \cong \Cl_{r,s+8} \cong M_{16} (\Cl_{r,s})\,.\]
\end{theorem}

\Proof: Since $\Cl_{1,1}$ is positive we get
\[\Cl_{r+1,s+1} \cong
\Cl_{r,s} \otimes \Cl_{1,1} \cong \Cl_{r,s} \otimes M_2(\R)\]
and
repeatedly by Theorem 3
\[\Cl_{r+n,s+n} \cong \Cl_{r,s} \otimes
M_2(\R) \oco_{n-{\rm mal}} M_2(\R) \cong \Cl_{r,s} \otimes
M_{2^n}(\R) \cong M_{2^n}(\Cl_{r,s})\,.\]
Since $\Cl_{4,0}$ is
positive again by Theorem 3 and by the Proposition we have
\[\Cl_{8,0} \cong \Cl_{4,0} \otimes \Cl_{4,0} \cong\Cl_{4,0} \otimes
\Cl_{0,4} \cong \Cl_{4,4} \cong M_{16}(\R)\]
and
\[\Cl_{p+8,q}
\cong \Cl_{p,q} \otimes \Cl_{8,0} \cong \Cl_{p,q} \otimes M_{16}(\R)\,,\] 
because $\Cl_{8,0}$ is positive, too.
\endproof 

We end up with the following table displaying the special Clifford algebras
$\Cl_n =\Cl_{n,0}$

\hfuzz=3pt

$\vbox{\tabskip=0pt \offinterlineskip
        \halign{\tabskip=2.5pt plus1em#%
                 &\strut\ctr{$#\quad$}
                &#\vl\quad&\ctr{$\,#\,$}
                &\ctr{$\,#\,$}
                &\ctr{$\,#\,$}
                &\ctr{$\,#\,$}
                &\ctr{$\,#\,$}
                &\ctr{$\,#\,$}
                &\ctr{$\,#\,$}
                &\ctr{$\,#\,$}
                &\ctr{$\,#\,$}
                          \cr
&n&&0 &1&2&3&4&5&6&7&8 \cr \noalign{\hrule}
&\Cl_n&&\R&\C&\H&\H\oplus \H&M_2(\H)&M_4(\C)&M_8(\R)&M_8(\R)\oplus
M_8(\R)& M_{16}(\R)\cr }}$
\bigskip
\bigskip

\section{\bf Representations of Clifford algebras}
\bigskip

We also need representations of abstract Clifford algebras. Recall
that a representation $\rho:\Cl_n\to \End(E)$ on a real (or
complex) finite dimensional vector space $E$ is irreducible if for
any decomposition $E=E_1\oplus E_2$ into subspaces invariant under
$\rho$ one has $E_1=E$ or $E_2=E$. In the reducible case one has
$\rho=\rho_1\oplus \rho_2$ with $\rho_j=\rho|_{E_j}$. Give any
non-trivial representation $\rho$, one can find an inner product
$\l\cdot,\cdot\r$ on $E$ such that $\rho(x)$ acts orthogonally (or
unitarily) on $E$ for all $x\in\R^n\subset\Cl_n$ with $|x|=1$. One
merely has to average a given inner product $\l\cdot,\cdot\r'$
over the finite (multiplicative) group $G_n$ generated by
$e_1,\dots,e_n\in\Cl_n$, i.e.\ one puts
\[\l v,w\r=\sum_{x\in
G_n}\l\rho(x)v,\rho(x)w\r',\quad v,w\in E.\]
Since
$\rho(x)^2=-|x|^2 I_E$, this amounts to
\[\l\rho(x)v,w\r = - \l v,\rho(x)w\r, \quad v,~w\in E,~ x\in\R^n,\]
i.e.\ $\rho(x)^\ast
=-\rho(x)$. If this holds we call $\rho$ a skew-adjoint
representation.

Now any representation $\rho$ can easily be decomposed into a
direct sum of irreducible ones: Choosing $v\in E$, $v\ne 0$, one
considers $E_v=\{\rho(x)v\mid x\in\Cl_n\}$ which is invariant
under $\rho$. Since $E_v^\bot$ is also invariant, successively
splitting off invariant subspaces (in case also of $E_v$) one ends
up with $E=\bigoplus_{j=1}^m E_j$ and
$\rho=\bigoplus_{j=1}^m\rho_j$ where $\rho_j$ is irreducible. Two
representations $\rho_j:\Cl_n\to \End(E_j)$ are called equivalent
if they are implemented by an isomorphism $T:E_1\to E_2$, i.e.\
\[T\rho_1(x)=\rho_2(x)T,\quad\hbox{for all}\quad x\in\Cl_n.\]
In our case we have $\Cl_n$ of the form $M_m(\K)$ if $n\ne 3$ and
7 which being a simple algebra does not contain any non-trivial
two-sided ideal. To see this consider elementary matrices $e_{ij}$
with entries 1 at $i,j$ and 0 elsewhere. Now given a two-sided
ideal $V\subset M_m(\K)$ and $x=\sum_{1\le i,j\le
m}x_{ij}e_{ij}\in V\setminus \{0\}$ there is an $x_{ij}\ne 0$ and
therefore $e_{ij}=x_{ij}^{-1}e_{ii}xe_{jj}\in V$. Since
$e_{ij}e_{k\ell}=\delta_{kj}e_{i\ell}$ all of the $e_{ij}$ belong
to $V$, i.e.\ $V=M_m(\K)$. In particular, we consider the left
regular representation
\[\rho_L:M_m(\K)\to\End\big(M_m(\K)\big)\quad \rho_L(x)y=xy,~
x,y\in M_m(\K).\] It decomposed as
$\rho_L=\bigoplus_{j=1}^m\rho_j$ with irreducible representations
$\rho_j(x)ye_{jj}=x y e_{jj}$ on the left ideals
$V_j=M_m(\K)e_{jj}$. If $\rho$ is an arbitrary faithful (i.e.\
injective) irreducible representation it has to be equivalent to
one of the $\rho_j$ and hence to each of them. To prove this note
that there is a $v\in E$ and an $x\in V_1$ with $\rho(x)v\ne 0$.
Now define $T:V_1\to E$ by $T(y)=\rho(y)v$, $y\in V_1$, and
observe that
\[T\rho_1(z)y=T(z y)=\rho(z
y)v=\rho(z)\rho(y)v=\rho(z)T y,\quad y\in V_1,\]
hence by Schur's
Lemma $T$ has to be an isomorphism since both representations are
irreducible: $\ker T\subset V_1$ and ${\rm im}~ T\subset E$ are
subspaces invariant under $\rho_1$ and $\rho$, respectively, hence
${\rm im}~ T=E$ and $\Ker T=\{0\}$, since $T\ne 0$. Combining this
with Theorem 4 and the table above we obtain:

\begin{theorem} For $n\not\equiv 3$ and 7 $\mod(8)$ the
Clifford algebra $\Cl_n$ has up to equivalence exactly one
irreducible representation, viz.\ on $\R^{a_n}$, where
\[a_n=\cases{1,\quad &$n=0$,\cr 2, &$n=1$,\cr 4, &$n=2,3$,\cr 8,
&$n=4,5,6,7$,\cr}\]
and $a_{n+8k}=2^{4k}~a_n$. In cases $n\equiv
3$ or 7 $\mod(8)$ there are exactly two non-equivalent irreducible
representations on $\R^{a_n}$.
\end{theorem}

\Proof: Noting that $\C$ is irreducibly represented in $M_2(\R)$ by 
$a+ib\mapsto \pmatrix{a&-b\cr b& a\cr}$ and $\H$ in $M_2(\C)\subset M_4(\R)$ 
by $z+wj\mapsto \pmatrix{z &w\cr -\bar w &\bar z\cr}$ the first assertion 
follows from the table. For $n=3$ or 7 one has $\Cl_n\cong M_n(\K)\oplus 
M_n(\K)$ with $\K=\R$ or $\H$, respectively, and two irreducible 
representations on $\K^n\simeq\R^{a_n}$ are given by $\rho_1(x,y)=\rho(x)$ 
and $\rho_2(x,y)=\rho(y)$. They are not equivalent since 
$\rho_1(I_n,-I_n)=I_n$ and $\rho_2(I_n,-I_n)=-I_n$. 
\endproof
\medskip

Writing $n=(2\ell+1)16^\alpha2^\beta$ with $\beta=0,1,2$, or $3$
and $\rho(n)=8\alpha+2^\beta$ the highest power of 2 dividing $n$
being just $a_{\rho(n)-1}$ we obtain:
\medskip

\noindent
{\bf Corollary}~ {\it The Clifford algebra $\Cl_{\rho(n)-1}$ has a
non-trivial representation on $\R^n$. In particular, there are
matrices $A_1,\dots,A_{\rho(n)-1}\in O(n)$ with $A_i^2=-I_n$ and
$A_iA_j=-A_jA_i$, $i\ne j$, $i,j=1,\dots,\rho(n)-1$.}
\smallskip

\Proof: Let $n=p\cdot a_{\rho(n)-1}$, $p$ odd, and
$\delta:\Cl_{\rho(n)-1}\to \End(\R^{a_{\rho(n)-1}})$ the previous
representation. Then
\[\bar\delta=\bigoplus_{k=1}^p\delta:\Cl_{\rho(n)-1}\to
\End\Big(\bigoplus_{k=1}^p\R^{a_{\rho(n)-1}}\Big)\]
is the one we
are looking for, since, as seen before, we can choose an inner
product that renders $A_i$ orthogonal with respect to a suitable
orthonormal basis.
\endproof
\medskip

The matrices $A_j$ and the numbers $a_{\rho(n)-1}$ which are guaranteed by 
the Corollary are often called Hurwitz-Radon matrices and Radon numbers, 
respectively, after A.\ Hurwitz [Hur] and J.\ Radon [Rad] who around 1920 
independently constructed such matrices in order to factorize quadratic 
forms. They also solved the linear vector field problem: There are exactly 
$a_{\rho(2n)-1}$ linear vector fields, given by $X_j(x)=A_jx$, $x\in 
S^{2n-1}\subset\R^{2n}$, that are linearly independent at each point; cf.\ 
[Eck]. 

We also consider complex Clifford algebras
$\Cl^\C_n=\Cl_n\otimes_\R\C$ and their irreducible representations
on complex vector spaces. Complexifying immediately entails
\[\Cl_n^{\C}\cong\cases{M_{2^{k}}(\C),\quad &if $n=2k$,\cr
M_{2^{k}}(\C)\oplus M_{2^{k}}(\C), & if $n=2k+1$.\cr}\]
This also
shows that (up to equivalence) $\Cl_n^\C$ has exactly one
irreducible representation if $n=2k$ and exactly two if $n=2k+1$.
The isomorphism with $M_{2^k}(\C)$ can be made explicit using the
Pauli matrices $\sigma_j$. The basis elements $e_j$, $1\le j\le
2k$, are represented (up to a choice of sign) by the following
skew-hermitian unitary matrices:
\begin{eqnarray*}
A_{2\ell-1} = &iA_{2\ell-1}'&=\sigma_3\moco_{\ell-1\hbox{-}times}
\sigma_3\otimes i\sigma_1\otimes I_2\moco_{n-\ell\hbox{-}times} I_2,~1\le
\ell\le k,\\
A_{2\ell} = &iA_{2\ell}'&=\sigma_3\moco_{\ell-1\hbox{-}times}
\sigma_3\otimes i\sigma_2\otimes I_2\moco_{n-\ell\hbox{-}times}
 I_2,~1\le \ell\le k.
\end{eqnarray*}
This is a simple consequence of the construction in Theorem 3.

These matrices allow to classify complex Clifford algebras and
their irreducible representations directly. If $n=2k$, i.e., $\dim
M_{2^k}(\C)=2^{2k}=\dim\Cl_n$ one only has to show that the
representation $\rho(e_j)=A_j\in M_{2^k}(\C)$, $j=1,\dots,2k=n$,
is faithful, i.e.\ that the matrices $A_I=A_{i_1}\cdots
A_{i_\ell}$, $1\le i_1<\cdots< i_\ell\le 2k$, are linearly
independent. To this end one uses the trace which defines an inner
product on $M_{2^k}(\C)$ by $\l A,B\r=\tr(A^\ast B)$. Now for
$\ell$ even one has
\[\tr(A_I)=\tr(A_{i_\ell}A_{i_1}\cdots
A_{i_{\ell-1}}=(-1)^{\ell-1}\tr(A_I),\] hence $\tr(A_I)=0$, and
for $\ell<2k$ odd and $i_{\ell+1}\not\in I$ one has
\[\tr(A_I)=-\tr(A_I A_{i_{\ell+1}}A_{i_{\ell+1}})
=-\tr(A_{i_{\ell+1}}A_IA_{i_{\ell+1}})=(-1)^{\ell}\tr(A_I),\]
hence again $\tr(A_I)=0$. Given a linear combination $\sum
a_IA_I=0$ this implies
\[0=\tr\Big(\sum a_IA_IA_J\Big)=\pm a_J 2^k.\]
Note that the argument does not use the special shape of the
matrices $A_j$.\\ If $n=2k+1$ there is another matrix
\[A_{2k+1}=-i A_{2k+1}'=-i\sigma_3\moco\sigma_3.\]
However, the extended representation $\rho:\Cl_{2k+1}\to
M_{2^k}(\C)$ by $\rho(e_{2k+1})=A_{2k+1}$ is no longer faithful,
since
\[\omega=i^{[(n+1)/2]}e_1\cdots e_n=i^{k+1}e_1\cdots
e_{2k+1}\]
is represented by $\rho(\omega)=I_{2^k}$.\\ A
non-equivalent representation $\rho'$ will be defined by
$\rho'(e_j)=-A_j$, $1\le j\le 2k+1$, since
$\rho'(\omega)=-I_{2^k}$. To obtain a faithful (reducible)
representation one takes the direct sum
$\rho\oplus\rho':\Cl_{2k+1}^\C\to M_{2^k}(\C)\oplus
M_{2^k}(\C)\subset M_{2^{k+1}}(\C)$.

\begin{definition}\quad If $\rho:\Cl_n^\C\to\End(E)$ is an
irreducible faithful representation, then the vector space
$E\cong\C^{2^k}$ is called a space of spinors; usually, it will be
denoted by $S_0$.
\end{definition}

{\bf Remarks} 1. Different realizations of $S_0$ will be given in the 
following examples. 

2. In $\Cl_{2k}^\C$ the element $\omega=i^ke_1\cdots e_{2k}$ satisfies 
$\omega^2=1$ and $\omega\cdot e_j=-e_j\omega$, $j=1,\dots,2k$, hence defines 
a $\Z_2$-grading on $S_0$, i.e.\ $S_0=S_0^0\oplus S_0^1$ where 
$S_0^j=\frac{1}{2}\big(1+(-1)^j\omega\big)S_0$, $j=0,1$, are the so-called 
spaces of half-spinors. 

The uniqueness of irreducible representations by complex $2^k\times 
2^k$-matrices that contain and generalize Pauli's matrices [Pau] has first 
been proved by P.\ Jordan and E.\ Wigner [JW] using group theoretical 
arguments (in connection with the quantum theory of many electron systems in 
1927). The shortest proof without any theory of real Clifford algebras can be 
found in H.\ Weyl's ``Group Theory and Quantum mechanics'' of 1931. He 
explicitly gives the matrices $A_j'$ and expresses by them all of the 
elementary matrices that generate the simple algebra $M_{2^k}(\C)$; cf.\ also 
[BW] and [Wey]. We give his construction in the following example. 
\medskip

{\bf Examples} 3. The reducible representation
$\rho_L:\Cl_n^\C\to\End(\Cl_n^\C)$  be decomposed into a sum of
irreducible ones if $n=2k$. In the first case one needs a minimal
left ideal $V$ to act on. Starting from an orthonormal basis
$\{e_1,\dots,e_{2k}\}$ of $\C^n$ one can construct $V$ as follows:
Put
\[f_\ell=\frac{1}{\sqrt{2}}(e_{2\ell-1}+i
e_{2\ell})\quad\hbox{and}\quad
g_\ell=\frac{1}{\sqrt{2}}(e_{2\ell-1}-i e_{2\ell})\] as well as
\[p_\ell^\pm=\frac{1}{2}(1\pm ie_{2\ell-1}e_{2\ell})\quad
\hbox{for}\quad 1\le\ell\le k,\]
hence
$p_\ell^+=-\frac{1}{2}f_\ell g_\ell$ and
$p_\ell^-=-\frac{1}{2}g_\ell f_\ell$. The idempotents $p_\ell^\pm$
mutually commute, and for any $n$-tuple $\e=(\e_1,\dots,\e_k)$
with $\e_j=\pm$ they define a projection $p^\e=p_1^{\e_1}\cdots
p_k^{\e_k}$. Then $V^\e=\Cl_{2k}^\C p^\e$ is a minimal left ideal
and $\Cl_{2k}^\C\cong \End(V^\e)$.

Note that each projection $p^\e$ is associated with an elementary matrix
$e_{jj}$ in $M_{2^k}(\C)$, e.g. $e_{11}$ with $p=p^\e$ where
$\e=(1,\dots,1)$. Thus $V$ is isomorphic with the vector space of matrices
that have non-trivial entries only in its $j^{th}$ column. Therefore, one has
$1=\sum_{\e}p^\e$.

With regards to this example B.L.\ van der Waerden writes in 1966 [vdW2]:
\begin{quote}
\baselineskip10pt {\footnotesize If you want to determine the structure of an 
algebra or of a group defined by generating elements and relations and to 
find a representation of the algebra or group by linear transformations or by 
permutations, construct the regular representation.} 
\end{quote}

4. To decompose the reducible representation $c:\Cl_n^{\C}\to \End(\L\C^n)$ 
which in fact is equivalent to the previous one one starts with the 
orthogonal decomposition $\C^n=W\oplus\o{W}$, where $W$ or $\o{W}$ denote the 
subspaces spanned by $g_\ell$ or $f_\ell$, respectively. From the relations 
\begin{eqnarray*}
&&f_jf_\ell+f_\ell f_j=0,\\ &&g_jg_\ell+g_\ell g_j=0,\\ &&f_jg_\ell+g_\ell
f_j=-2\delta_{j\ell}.
\end{eqnarray*}
and since $f_jg_Ip=0$ if $j\not\in I$ and $f_j
g_Ip=(-1)^{\ell}2g_{I'}p$ if $I=I'\cup\{j=i_\ell\}$ one obtains
that the subspace
\[V=\Cl_{2k}^\C p=\L W p\]
is a left ideal and isomorphic with $\L W$ as a vector space.
Modifying $c$ on $\L W\subset\L\C^{2k}$ by taking
\[\tilde
c(w)=\sqrt{2}\big(\e(v)-i(\bar v)\big)\in\End(\L W)\] for
$w=v+\bar v\in W\oplus \o{W}$ one obtains the appropriate
irreducible representation. Indeed, from the previous relations
one easily verifies for
\[w=\sum_{j=1}^k(x_je_{2j-1}+y_je_{2j})
=\frac{1}{\sqrt{2}}\sum_{j=1}^k(\bar{z_j}f_j+z_jg_j)\]
with
$z_j=x_j+iy_j$ the relation $\tilde c(w)^2=-\sum_{j=1}^k|z_j|^2
I=-|w|^2 I$.

The main problem with the space of spinors is that there is no canonical way 
to decompose a given representation, even a natural one as in the previous 
examples, into irreducible ones. Therefore, the spin structure to be defined 
in the next section and whose construction rests on a proper choice of 
irreducible representations will in its last analysis always be superficial. 

We have shown that up to equivalence any representation
$\rho:\Cl_{2k}\to\End(E)$ can be written as $\rho_0\otimes
I:\Cl_{2k}\to\End(S_0\otimes W)$ with $E\cong S_0\otimes W$ and
where
\[\rho(v)(e\otimes w)=\rho_0(v)e\otimes w,\quad e\otimes
w\in S_0\otimes W.\]
Now given $S_0$, at least, $W$ is canonically
defined. To see this we have to digress and recall some general
results about tensor products.

Given two real (or complex) vector spaces $E$ and $F$, which moreover are 
right respectively left modules for some real (or complex) algebra $A$ the 
tensor product $E\otimes_A F$ is defined as the quotient space of $E\otimes 
F$ by the subspace generated by $va\otimes w-v\otimes a w$, $v\in E$, $w\in 
F$, $a\in A$. It is the unique vector space with the following universal 
property. If $H$ is a vector space and $f:E\times F\to H$ is a bilinear 
$A$-balanced map, i.e. $f(v a,w)=f(v,a w)$ for $v\in E$, $w\in F$, $a\in A$, 
then there is a unique linear map $f_A$ such that the following diagram 
commutes: 

\unitlength1cm 
\begin{picture}(6,4)(-4,-1)
\put(-0.4,0){$E\times F$} \put(0.7,0.1){\vector(1,0){3.1}} \put(4,0){$H$}
\put(0.1,0.5){\vector(0,1){1.5}} \put(-0.5,2.3){$E\otimes_A F$}
\put(0.5,2){\vector(2,-1){3.3}}
\put(-0.2,1.3){$\scriptstyle{\gamma}$}
\put(1.7,0.4){$\scriptstyle{f}$}
\put(1.7,1.7){$\scriptstyle{f_A}$}
\end{picture}

Let $B$ denote another algebra and let $G$ be a left-$B$-module.
Then the following results hold:\\ (a) If $F$ is also a
right-$B$-module (hence an $(A,B)$-bimodule), then $E\otimes_A F$
is a right-$B$-module, $F\otimes_B G$ an left-$A$-module, and
\[(E\otimes_A F)\otimes_B G\cong E\otimes_A(F\otimes_B G)\]
is a natural isomorphism. \\(b) If $E$ is a $(B,A)$-bimodule, then
$\Hom_B(E,G)$, the space of $B$-module homomorphisms consisting of
linear maps $f\in\Hom(E,G)$, which satisfy $f(b v)=b f(v)$ for
$v\in E$ and $b\in B$, is an left-$A$-module by $a f(v)=f(v a)$,
$a\in A$, $v\in E$. One has the natural isomorphism
\[\Hom_A\big((F,\Hom_B(E,G)\big)\cong \Hom_B(E\otimes_A F,G),\]
induced by $f\mapsto\tilde f$ mit $\tilde f(u\otimes v)=f(v)(u)$
for $f\in \Hom_A(F,\Hom_B(E,G))$, $u\in E$, and $v\in F$.\\(c)
Moreover, by $(f\otimes v)(w)=f(w)\otimes v$ for
$f\in\Hom_B(G,E)$, $v\in F$, and $w\in G$, one obtains a natural
homomorphism
\[\Hom_B(G,E)\otimes_A F\cong \Hom_B(G,E\otimes_A F),\]
which is one-to-one and onto if $F$ is a finitely generated
projective module, i.e.\ a direct summand of $A^n$ for some $n\in
\N$.

(d) If, however, $E$ is an right-$A$-module, $G$ a
left-$B$-module, and $F$ a $(B,A)$-bimodule, then $\Hom_A(E,F)$ is
a left-$B$-module by $(b f)(v)=b f(v)$, $b\in B$, $v\in E$, and
one has a natural isomorphism
\[E\otimes_A\Hom_B(F,G)\cong \Hom_B\big(\Hom_A(E,F),G\big),\]
induced by $u\otimes f\mapsto
\tilde h$ with $h(g)=f\circ g(u)$ for $u\in E$, $f\in
\Hom_B(F,G)$, and $g\in \Hom_A(E,F)$.\\We only need these results
in the special case $A=\C$ and leave its proofs to the reader;
cf.\ [AF].

As a simple consequence of the last one we obtain that the module
$W$ in the decomposition $E=S_0\otimes W$ can be chosen as
$W=\Hom_{\Cl_n^\C}(S_0,E)$: Since $\Cl_n^\C=\End(S_0)$, there are
isomorphisms
\[S\otimes\Hom_{\Cl_n^\C}(S_0,E)\cong
\Hom_{\Cl_n^\C}\big(\Hom(S_0,S_0),E\big)\cong
\Hom_{\Cl_n^\C}(\Cl_n^\C,E)\cong E.\] 
We conclude this section and the 
algebraic part of the paper with a classical result of representation that 
will be essential in the proof of the main theorems of the next section. It 
is a special case of the Theorem of Skolem-Noether. 
\medskip

\noindent 
{\bf Lemma}\quad {\it Let $A_1$ and $A_2$ be two isomorphic simple 
subalgebras of $M_n(\C)$, say both isomorphic to $M_k(\C)$. Then each 
isomorphism $\Phi:A_1\to A_2$ is an inner automorphism ${\rm Ad}(U)$ of 
$M_n(\C)$, i.e., there is a $U\in GL_n(\C)$ with 
\[\Phi(a)={\rm Ad}(U)a=UaU^{-1},\quad a\in A_2.\]
In particular, each automorphism of $M_k(\C)$ is inner and each
derivation $D$ of $M_k(\C)$ is an inner derivation, i.e.\ given by
\[D a={\rm ad}(v)a=[v,a]=v a-a v,\quad a\in M_k(\C),\]
for some $v\in M_k(\C)$.}  
\smallskip

\Proof: The simple algebra $M_k(\C)$ is represented by $A_1$ and
$A_2$ in $M_n(\C)$, respectively. There are decompositions
$\C^n=\bigoplus_{j=1}^\ell E_j$ and $\C^n=\bigoplus_{j=1}^r F_j$
which reduce $A_1$ and $A_2$, respectively. $A_1$ and $A_2$ being
isomorphic, one has $\ell = r$, and since the restricted
irreducible representations have to equivalent, one has
\[E_j\cong \C^k\cong F_j\quad\hbox{for all}~j.\]
Now $U$ is given as the
direct sum of such isomorphisms. To prove the second assertion,
one simply has to take $k=n$ and to choose $A_2$ as the image of
$A_1=M_k(\C)$ under a given automorphism. For the last assertion
take $n=2k$, $A_1=\left\{\pmatrix{a&0\cr 0&a\cr}\mid a\in
M_k(\C)\right\},$ and $A_2=\left\{\pmatrix{a&D(a)\cr 0&a\cr}\mid
a\in M_k(\C)\right\}$ for a given derivation $D$. Then there is a
$U=\pmatrix{u&v\cr w&z\cr}$ with
\[\pmatrix{u&v\cr
w&z\cr}\pmatrix{x&0\cr 0&x\cr}= \pmatrix{x&D(x)\cr
0&x\cr}\pmatrix{u&v\cr w&z\cr}\]
for all $x\in M_k(\C)$. This
entails $wx=xw$ and $zx=xz$ for all $x\in M_k(\C)$, hence, by
Schur's Lemma, $w$ and $z$ are multiples of the identity. If say
$z\ne 0$, the further condition $xv+D(x)z=vx$ on $D$ leads to
$D={\rm ad}(z^{-1}v)$.
\endproof
\medskip

\noindent 
{\bf Remark}\quad The maps $\nu:GL_k(\C)\to \Aut(M_k(\C))$, 
$\tau(u)={\rm Ad}(u)$ and $\mu:M_k(\C)\to {\rm Der}(M_k(\C))$, $\mu(v)={\rm 
ad}(v)$ (into the space of derivations) are both onto but in general not 
one-to-one, since $\ker\nu\cong\C^\ast=\C\setminus\{0\}$ and $\ker 
\mu\cong\C$. In the case of $\mu$ one can, however, consider its restriction 
$\mu_0$ to the subspace $M_k^0(\C)$ of matrices with vanishing trace and 
obtains an isomorphism. 
\bigskip
\bigskip
\section{\bf Spinor bundles and Dirac operators}
\bigskip
We now want to globalize the results of the previous section, i.e.\ to 
perform the constructions on vector bundles over smooth manifold. 
\medskip

\begin{definition} Let $E$ be a Euclidean vector bundle of
rank $k$ over $M$. The vector bundle $\Cl(E) = \coprod_{p\in M}
\Cl(E_p)$ will be called the Clifford bundle of $E$. If $M$ is
endowed with a Riemannian structure one particularly has $\Cl M=
\Cl(TM)$, the Clifford bundle of $M$.
\end{definition}

Starting from a local orthonormal frame $(\fe_i)_{1\le i\le k}$ of
$E$ over $U$ one obtains a local trivialization
\[E|_U\ni
v_q=\sum_{i=1}^k a_i\fe_i(q)\mapsto \vi(v_q)=\big(q,(a_i)_{1\le
i\le k}\big)\in U\times\R^k\] and
$\vi|_{E_q}:E_q\to\{q\}\times\R^k=\R^k$ is isometric for any $q\in
U$. Choosing an atlas $\ca=\{(U_\alpha,\vi_\alpha)\mid\alpha\in
A\}$ in this way one obtains a cocycle of transition maps with
$g_{\alpha\beta}:U_\alpha\cap U_\beta\to O(k)$. To trivialize
$\Cl(E|_U)$ we choose
\[\Cl(\vi_\alpha)(x_q)=\Cl(\vi_\alpha|_{E_q})(x_q)\in\{q\}\times\Cl_k,\quad
x_q\in\Cl(E_q)\] according to Remark 1 after Theorem 1. The
corresponding transition maps are given by
$f_{\alpha\beta}:U_\alpha\cap U_\beta\to \Aut(\Cl_k)$ with
\[f_{\alpha\beta}(p)=\Cl\big(g_{\alpha\beta}(p)\big),\quad p\in
U_\alpha\cap U_\beta.\] They possess the cocycle property and are
differentiable, since the group homomorphism $O(k)\ni f\mapsto
\Cl(f)\in \Aut(\Cl_k)$ is a polynomial in the coefficients with
respect to a fixed orthonormal basis of $\R^k$ and the induced
basis of $\Cl_k$. This makes $\Cl(E)$ a smooth vector bundle. In
particular,
\[\{\fe_{i_1}\cdot \cdots \cdot \fe_{i_k}\mid 1\le i_1
< \cdots < i_k \le r, \; k=0,\cdots,n\}\]
provides an orthonormal
frame of $\Cl(E)\big|_U$. The Clifford bundle $\Cl(E)$ depends on
the Euclidean structure and is itself a Euclidean vector bundle.
On the other hand, the $C^\infty$-structure and the Riemannian
structure of $\Cl(E)$ do not depend on the choice of the frame.

Each fiber of $\Cl(E)$ comes with an algebra structure and
fiber-wise multiplication makes $C^\infty\big(\Cl(E)\big)$, the
space of smooth sections, into an algebra, too. Suitably modifying
the definition of a vector bundle one obtains the notion of an
algebra bundle $(A,\pi,M)$:\\Each fiber $\pi^{-1}(p)$ is a finite
dimensional topological algebra with respect to the topology
induced by $A$, and at each point $p\in M$ there exists a chart
$\vi:\pi^{-1}(U)\to U\times A$ with a fixed given algebra $A_0$,
such that
\[\vi|_{\pi^{-1}(q)}:\pi^{-1}(q)\to\{q\}\times A_0\]
is an algebra isomorphism for any $q\in U$. We only consider the
special case of unital algebras $A_0$ and $A_p$ with unit elements
$e_0$ and $e_p$, respectively. Then we have a global section $\fe$
in $A$.\\Alternatively, we may assume a bundle morphism
$(\mu,\id_M)$ with $\mu:A\otimes A\to A$ and $\mu(e_p\otimes
a_p)=a_p=\mu(a_p\otimes e_p)$ for any $a_p\in A_p$.\\Moreover a
vector bundle $F$ will be called a (left-)$A$-bundle if there is a
bundle morphism $\tau:A\otimes F\to F$ with
\[\tau\big(a_p\otimes\tau(b_p\otimes
v_p)\big)=\tau\big(\mu(a_p\otimes b_p)\otimes
v_p\big)\quad\hbox{for}\quad a_p,b_p\in A_p,~ v_p\in F_p.\] In
other words, $F_p$ is a left-$A_p$-module for any $p\in M$ and
\[\sigma\cdot \fs(p) = \sigma(p) \cdot \fs(p) = \tau\big(\sigma(p)
\otimes \fs(p)\big)\,,\quad p\in M\]
defines a smooth section,
$\sigma\cdot \fs\in C^\infty(F)$, for $\sigma \in C^\infty(A)$ and
$\fs\in C^\infty(F)$. An $A$-bundle morphism is a bundle morphism
$(\hat f,\id_M)$ between two $A$-modules $F$ and $G$, that is an
$A_p$-linear map from $F_p$ to $G_p$ for any $p\in M$. The space
of $A$-bundle morphisms, $\HOM_A(F,G)$, is a $C^\infty(A)$-module
and can be identified with $C^\infty\big(\Hom_A(F,G)\big)$. Here
$\Hom_A(F,G)$ is the sub-bundle of $\Hom(F,G)$, whose fibers are
$\Hom_{A_p}(F_p,G_p)$, $p\in M$. In particular, we can extend the
natural isomorphisms at the end of the previous section to the
setting of $A$-bundles.

\noindent We are now going to define geometric differential operators that
are closely connected with the topological or geometrical structure of an
oriented Riemannian manifold $M$.

\begin{definition}\quad A smooth vector bundle $E$ over $M$ is called a
spinor bundle over $M$ if it is a left-$\Cl M$-bundle.
\end{definition}

If the module structure is given by the morphism $\tau : \Cl M
\otimes E\to E$ we also consider the bundle morphism
$c_E:TM\to\End(E)$, $c_E(v_p)(e_p)=\tau(v\otimes e)$, $v\in T_pM$,
$e_p\in E_p$, induced by $\tau$ and its extension $c_E:\Cl M\to
\End(E)$ to a morphism of algebra bundles. To emphasize the
underlying Clifford multiplication we sometimes denote a spinor
bundle by $(E,c_E)$.
\medskip

{\bf{Examples}}\quad 5. The Clifford bundle $\Cl M$ itself is a
spinor bundle if
\[c_{\Cl M} : \Cl M \to \End(\Cl M)\]
is in each
fiber given by the left regular representation $\rho_L$.

6. Likewise the Grassmann bundle $\L^\ast M$, the exterior bundle
of the cotangent bundle $T^\ast M$, is turned into a spinor bundle
using the isomorphism of $\Cl M$ with $\L^\ast M$. Here and in the
following we use the ``musical isomorphisms'' $^\sharp:\L^\ast
M\to\L M$ and its inverse $^\flat:\L M\to\L^\ast M$ that extend
the pairing between tangent vectors and cotangent vectors provided
by the Riemannian metric $\fg$ of $M$.

7. Given a spinor bundle $E$ and a smooth vector bundle $F$ we can turn
$E\otimes F$ into a spinor bundle, $E$ twisted by $F$. Here $\Cl M$ operates
on $E\otimes F$ by $v\cdot(e\otimes f)=(v\cdot e)\otimes f$ for $e\otimes
f\in E\otimes F$.
\medskip

Given a spinor bundle $E$ over $M$ via the isomorphism $^\sharp:
T^\ast M\to TM\subset\Cl M$ the bundle morphism $\tau$ induces a
linear map $T :C^\infty(T^\ast M \otimes E)\to C^\infty(E)$ given
by
\[T(\omega\otimes \fs)(p) = \tau\big(\omega(p)^\sharp \otimes
\fs(p)\big),~ p\in M.\]
It is easy to see that $T$ is a
differential operator of order zero.
\medskip
To define more sophisticated differential operators on
$C^\infty(E)$ we need a (Koszul) connection $\nabla$ on $E$, i.e.\
a linear differential operator $\nabla : C^\infty(E) \to
C^\infty(T^\ast M \otimes E)$ satisfying
\[\nabla (f\fs) = df
\otimes \fs + f\nabla \fs,~f\in C^\infty(M),~
 \fs\in C^\infty(E).\hfill(\ast)\]
For a vector field $X\in C^\infty(TM)$ this gives rise to a
covariant derivative $\nabla_X$ that satisfies
\[\nabla_X(f\fs)=X(f)\fs+f\nabla_X(\fs).\]
The dual connection
$\nabla^\ast$ on $C^\infty(E^\ast)$ can be defined by its
covariant derivatives
\[\nabla^\ast_X \fs^\ast(\fs) = X\big(\fs^\ast(\fs)\big) -
\fs^\ast(\nabla_X\fs)\] for $\fs^\ast\in C^\infty(E^\ast)$,
$\fs\in C^\infty(E)$, and $X\in C^\infty(TM)$. We also note the
following elementary constructions that can be performed with
connections $\nabla^E$ and $\nabla^F$ for vector bundles $E$ and
$F$, respectively. By
\[\nabla^{E\oplus F} (\fs\oplus \ft) =
\nabla^E \fs \oplus \nabla^F \ft\,,\]
and by
\[\nabla^{E\otimes
F}(\fs\otimes \ft) = (\nabla^E \fs) \otimes \ft + \Psi\big(\fs\otimes 
(\nabla^F\ft)\big)\,,\quad \fs\in C^\infty(E),\; \ft\in C^\infty(F)\] one 
defines connections $\nabla^{E\oplus F}$ for $E\oplus F$ and 
$\nabla^{E\otimes F}$ for $E\otimes F$. Here $\Psi$ is induced by the 
isomorphism of vector bundles, $\psi : E\otimes T^\ast M \otimes F \to T^\ast 
M \otimes E \otimes F$. In particular, one obtains a connection 
$\nabla^{\End(E)}$ on $\End(E)\cong E^\ast\otimes E$. 

\begin{definition} Let $E$ be a spinor bundle over $M$, and $\nabla$
a connection for $E$. Then
\[D= T\circ \nabla : C^\infty(E) \to C^\infty(E)\]
defines a first order differential operator, the
Dirac operator associated with $(E,\nabla)$.
\end{definition}

\begin{propo} Given a local orthonormal frame
$(E_i)_{1\le i\le m}$ of $TM$ over $U$ one has
\[D\fs = \sum_{k=1}^m E_k \cdot \nabla_{E_k} \fs\,\]
for $\fs\in
C^\infty(E\big|_U)$.
\end{propo}

\Proof: Since $X=\sum_{k=1}^m\l E_x,X\r E_k$ for $X\in C^\infty(TM)$ one has 
\[\nabla_X\fs =\sum_{k=1}^m\l E_k,X\r \nabla_{E_k}\fs,\]
hence
\[\nabla\fs=\sum_{k=1}^m E_k^{\flat}\otimes
\nabla_{E_k}\fs\]
and so the representation of $D$ as stated.
\endproof
\medskip

With respect to a local frame $(s_j)_{1\le j\le r}$ of $E$ a
connection is given by
\[\nabla\Big(\sum_{j=1}^r f_j\fs_j\Big) =
\sum_{j=1}^r \Big(d f_j\otimes \fs_j + f_j \sum_{k=1}^r
\omega_{jk} \otimes \fs_k\Big),\]
where the local connection form
$\omega=(\omega_{jk})_{1\le j,k\le r}$ defined on say $U$ uniquely
determines $\nabla$ on $U$ and vice versa.
\medskip

Recall that the tangent bundle of a Riemannian manifold $M$ itself
comes with a unique torsion-free Riemannian connection, the
Levi-Civita connection which we denote by $\o{\nabla}$. Here
torsion-free means that
\[\o{\nabla}\nolimits_X Y - \o{\nabla}\nolimits_Y X = [X,Y]\]
and Riemannian that
\[\l \o{\nabla}\nolimits_X Y,Z \r  + \l
Y,\o{\nabla}\nolimits_XZ\r  = X\l  Y,Z\r\] for any vector fields
$X$, $Y$ and $Z$. Moreover, the Levi-Civita connection
$\o{\nabla}$ extends to $T^\ast M$ and to the tensor bundle by the
previously mention constructions and also to the exterior bundle
$\L^\ast M$ and to the Clifford bundle if we assume the product
formula
\[\o{\nabla}\nolimits_X(\omega_1\wedge \omega_2) =
(\o{\nabla}\nolimits_X \omega_1)\wedge\omega_2
+\omega_1\wedge(\o{\nabla}\nolimits_X\omega_2)\] for forms
$\omega_1$, $\omega_2\in \mit\Omega(M)$ respectively
\[\o{\nabla}\nolimits_X(\sigma_1\cdot \sigma_2) =
(\o{\nabla}\nolimits_X\sigma_1)\cdot\sigma_2+
\sigma_1\cdot(\o{\nabla}\nolimits_X\sigma_2)\] for sections
$\sigma_1$, $\sigma_2\in\C^\infty(\Cl M)$. Combined with the
action of the Clifford bundle we obtain Dirac operators that are
defined on any oriented Riemannian manifold. The Dirac operator on
$\mit\Omega(M)$ has been introduced by E.\ K\"{a}hler in 1961 [Kae]
and so is sometimes called Dirac-K\"{a}hler. The extension
$\o{\nabla}$ to $\L^\ast M$ also satisfies
\[\o{\nabla}(\sigma\cdot \omega) = (\o{\nabla} \sigma) \cdot
\omega + \sigma \cdot (\o{\nabla} \omega)\] in the sense that
\[\o{\nabla}\nolimits_X(\sigma\cdot \omega) = (\o{\nabla}\nolimits_X
\sigma) \cdot \omega + \sigma \cdot (\o{\nabla}\nolimits_X
\omega)\]
for $X\in C^\infty(TM)$, $\omega \in \Omega(M)$, and
$\sigma\in C^\infty(\Cl M)$, hence in both cases $\o{\nabla}$ and
Clifford multiplication are compatible. Also recall that Clifford
multiplication by unit tangent vectors $X_p \in T_pM$ is
orthogonal on the spinor bundles $\Cl M$ and $\L^\ast M$ equipped
with the Riemannian metric induced by $\fg$. This suggests the
following definition.

\begin{definition}\quad Let $E$ be complex vector bundle with
a Hermitian metric $\l\cdot,\cdot\r$, a connection $\nabla$ and a
left $\Cl M_\C$-module structure $c_E$. We call the triple
$(E,\nabla,\l\cdot,\cdot\r)$ a Dirac triple and, for short, $E$ a
Dirac bundle if the given data are compatible, i.e.\ if

{\rm (1)} $c_E$ is a skew-adjoint representation in each fiber,

{\rm (2)} $\nabla$ is a compatible connection, i.e.
\[\nabla
(\sigma\cdot \fs) = (\o{\nabla} \sigma ) \cdot \fs + \sigma \cdot
\nabla \fs \quad\sigma\in C^\infty(\Cl M),~ \fs\in C^\infty(E),\]

{\rm (3)} $\nabla$ is a Riemannian connection, i.e.
\[\l\nabla_X
\fs_1,\fs_2 \r  + \l \fs_1,\nabla_X\fs_2\r  = X\big(\l
\fs_1,\fs_2\r\big),\quad X\in C^\infty(TM),~\fs_1,\fs_2\in
C^\infty(E).\]
\end{definition}

{\bf{Remarks}}\quad 1. By definition of the Levi-Civita connection
on $\Cl M$ to ensure (2) it suffices that
\[\nabla(X\cdot\fs)=(\o{\nabla} X)\cdot\fs + X\cdot\nabla\fs\]
for $X\in C^\infty(TM)\subset C^\infty(\Cl M)$ and $\fs\in
C^\infty(E)$.

2. If $(E,\nabla^E)$ is a Dirac bundle and $F$ is a Riemannian
vector bundle with Riemannian connection $\nabla^F$, then
$(E\otimes F,\nabla^E \otimes\nabla^F)$ with Clifford
multiplication as in Example 3 is again a Dirac bundle, since for
$\fs_1\in C^\infty(E)$, $\fs_2\in C^\infty(F)$, and $\sigma\in
C^\infty(\Cl M)$ one has
\begin{eqnarray*}
\nabla^E\otimes\nabla^F \big(\sigma\cdot (\fs_1\otimes\fs_2)\big) & = &
\nabla^E(\sigma\cdot\fs_1)\otimes\fs_2+(\sigma\cdot\fs_1)\otimes\nabla^F\fs_2\\
& = &\big((\nabla \sigma ) \cdot \fs_1\big)\otimes\fs_2 + (\sigma \cdot
\nabla^E \fs_1)\otimes\fs_2 + (\sigma\cdot\fs_1)\otimes\nabla^F\fs_2\\ & = &
(\nabla \sigma ) \cdot (\fs_1\otimes\fs_2) + \sigma \cdot (\nabla^E
\fs_1\otimes\fs_2) + \sigma\cdot(\fs_1\otimes\nabla^F\fs_2)\\ & = & (\nabla
\sigma ) \cdot (\fs_1\otimes\fs_2) + \sigma \cdot \nabla^E\otimes\nabla^F
(\fs_1\otimes\fs_2).
\end{eqnarray*}
Condition (1) also holds, since for $X\in C^\infty(TM)$
\begin{eqnarray*}
\big\l X\cdot (\fs_1\otimes \ft_1), \fs_2\otimes \ft_2 \big\r &=& \big\l
(X\cdot \fs_1)\otimes \ft_1, \fs_2\otimes \ft_2 \big\r \\ &=& \l X\cdot
\fs_1,\fs_2\r\l \ft_1, \ft_2 \r = -\l \fs_1,X\cdot \fs_2\r\l \ft_1, \ft_2 \r
\\ &=& -\big\l\fs_1\otimes \ft_1,(X\cdot\fs_2)\otimes\ft_2\big\r\\ &=& -\big\l
\fs_1\otimes \ft_1, X\cdot(\fs_2\otimes \ft_2)\big\r.
\end{eqnarray*}
In this way we obtain a Dirac operator with coefficients in the
bundle $F$ or a Dirac operator by twisting the Dirac operator
$D^E$ on $E$ with the connection $\nabla^F$. It will be denoted by
$D^E\otimes \nabla^F$ or simply by $D^E\otimes I_F$.
\medskip

It is well known that any complex vector bundle can be equipped with a 
Hermitian structure and with a Riemannian connection. Recall that one defines 
inner products and connections locally and in a second step uses partitions 
of unity to paste the local data to obtain global ones. So, in general, there 
is a lot of freedom to do this. In case of a complex spinor bundle one can 
ask whether these data can be chosen to satisfy (1) to (3). We shall prove 
that this can indeed be achieved. But before doing so we address the question 
of uniqueness, i.e.\ the impact that irreducibility has on the choice of 
these data. 
 
\begin{propo} Let $S$ be an irreducible complex
spinor bundle with a Hermitian metric $\l\cdot,\cdot\r$ and a
connection $\nabla$ satisfying properties (1) to (3). Then the
following results hold:\\ {\rm (a)} Any Hermitian metric
$\l\cdot,\cdot\r'$ with property (1) is of the form
\[\l\cdot,\cdot\r'=\lambda\l\cdot,\cdot\r\]
for some positive real-valued function $\lambda\in C^\infty(M)$.\\
{\rm (b)} Any connection $\nabla'$ with property (2) is of the
form
\[\nabla'=\nabla +\omega\]
for some complex-valued one-form
$\omega\in\Omega^1(M,\C)$.\\
{\rm (c)} If moreover $\nabla'$ is a
Riemannian connection with respect to the given metric, the
one-form $\omega$ is purely imaginary, i.e.\
$\nabla'=\nabla+i\eta$ for some real-valued one-form
$\eta\in\Omega^1(M,\R)$.
\end{propo}

\Proof: (a) For $p\in M$ let $T\in\End(S_p)$ be a hermitian endomorphism, 
such that 
\[\l s_1,s_2\r'=\l Ts_1,s_2\r\]
for all $s_1$, $s_2\in S_p$. Then
\[\l TX_p\cdot s_1,s_2\r=\l X_p\cdot s_1, s_2\r'
=-\l s_1,X\cdot s_2\r' =-\l Ts_1,X_p\cdot s_2\r=\l X_p\cdot
Ts_1,s_2\r\] for all $X_p\in T_pM$. Since the $X_p$ generate
$\End(S_p)$, $T$ commutes with each element of $\End(S_p)$, hence
by Schur's Lemma $T=\lambda I$ with $\lambda\in\C$. Since $T$ is
hermitian and positive, we have $\lambda\in\R$.\\(b) Analogously
we conclude that the section $\phi=\nabla'_X-\nabla_X$ into
$\End(S)$, which satisfies
\[\phi(\sigma\cdot\fs)=\sigma\cdot\phi(\fs)\]
for all $\sigma\in
C^\infty(\Cl M_\C)$ and $\fs\in C^\infty(S)$ because of the
derivation property that $\phi=\omega(X)I$ with
$\omega(X)\in\C$.\\(c) This is immediate, since
$\omega=\nabla'-\nabla$ has to be skew-hermitian, i.e.\
$\o{\omega}=-\omega$.
\endproof

\begin{theorem} Let $E$ be a complex spinor bundle over
the Riemannian manifold $M$ (of dimension $m=2n$). Then there are
a Hermitian structure and a Riemannian connection for $E$
compatible with Clifford multiplication which possess properties
(1) and (2).
\end{theorem}

\Proof: It suffices to prove this locally. Using a partition of unity local 
metrics as well as local connections can be pasted to global ones ensuing 
properties (1) to (3). Let $(U,\vi)$ be a chart of $M$ at $p\in M$ 
trivializing $E|_U$. We shall show that on a possibly smaller $U$ there are 
complex vector bundles $S$ and $W$ with $E|_U=S\otimes W$ and the $\Cl 
M$-action irreducible on $S$ and trivial on $W$. By the previous remarks it 
suffices to consider only $S$ and to equip $W$ with an arbitrary Hermitian 
structure and an arbitrary Riemannian connection.  Starting from a local 
orthonormal frame $\{E_1,\dots,E_m\}$ of $TM|_U$ we obtain sections 
$\fp^\e\in C^\infty(\Cl M_\C|_U)$ consisting of orthogonal projections. If 
$\fs_1\in C^\infty(E|_U)$ is a non-vanishing section one has 
$\fp^\e(q)\cdot\fs_1(q)\ne 0$ in the possibly smaller open set $U$ for some 
$\e$. Then 
\[f\big(q,(a_\sigma)_{\sigma\in G_m}\big) =\sum_{\sigma\in
G_m}a_\sigma \sigma(q)\fp^\e(q)\cdot\fs_1(q),\quad q\in U,~
(a_\sigma)_{\sigma\in G_m}\in\C^{|G|},\] defines a vector bundle
morphism $f:U\times \C^{|G_m|}\to E|_U$ of constant rank ${\rm rk}
f(p,\cdot)=N$ hence $F_1=\im f$ is a subbundle of $E$ whose fibers
are irreducible $\Cl_m^\C$-modules. We have $E=F_1\oplus F_1^\bot$
and proceeding likewise with a second non-vanishing section
$\fs_2\in C^\infty(F_1^\bot|_U)$ etc.\ we eventually obtain that
$E|_U\cong S\otimes W$ as a $\Cl U_\C$-bundle where $S=F_1$ and
$W=\e_U^\ell$. Now the products of sections $E_j$ in $C^\infty(\Cl
M|_U)$ generate a finite group $G_m$. Given an arbitrary Hermitian
structure $\l,\cdot,\cdot\r'$ on $E|_U$ we may define a new one by
putting
\[\l v,w\r=\sum_{\sigma\in G_m}\big\l\sigma(q)\cdot
v,\sigma(q)\cdot w\big\r',\quad v,w \in E_q.\]

Now given the irreducible spinor bundle $S$ we have an isomorphism
of algebra bundles $\Phi:\Cl M_\C|_U\to \End(S)$. Extending the
Levi-Civita connection $\o{\nabla}$ to $\Cl M|_U$ and then to $\Cl
M_\C|_U$, by $\Phi^{-1}$ we induce a connection
$\nabla=\Phi\o{\nabla}\Phi^{-1}$ on $\End(S)$. We only have to
show that $\nabla=\nabla^{\End(S)}$, i.e.\ induced by a Riemannian
connection $\nabla^S$ on $S$. This one will automatically possess
property (2), since
\begin{eqnarray*}
\nabla^S(\sigma\cdot\fs)&=&\nabla^S\big(\Phi(\sigma)(\fs)\big)
=\nabla\big(\Phi(\sigma)\big)(\fs)+\Phi(\sigma)(\nabla^S\fs)\\
&=&\Phi(\o{\nabla}\sigma)(\fs)+\Phi(\sigma)(\nabla^S\fs)\\
&=&\o{\nabla}\sigma\cdot\fs+\sigma\cdot\nabla^S\fs.
\end{eqnarray*}
Note that from the Remark concluding section 2 we have sub-bundles
$\End_0(S)$ of fiber-wise endomorphisms with trace 0 and ${\rm
Der}(S)$ of fiber-wise derivations of $\End(S)$, as well as a
bundle isomorphism $\mu_0:\End_0(S)\to{\rm Der}(S)$. \\ If
$\nabla_0$ is an arbitrary Riemannian connection on $S$ and
$\widetilde\nabla_0$ the connection induced on $\End(S)$, then
$\eta=\nabla-{\widetilde\nabla}_0$ is a section in $T^\ast
M\otimes\End\big(\End(S)\big)$ and from the derivation property
even a section in $T^\ast M\otimes{\rm Der}(S)$. For
$\gamma=\mu_0^{-1}\eta$ we then have
\[\nabla_X\ft-{\widetilde\nabla}_{0X}\ft=\gamma(X)\ft-\ft\gamma(X),\quad
X\in C^\infty(TM),~\ft\in C^\infty\big(\End(S)\big).\]
And putting
\[\nabla^S_{0X}\fs=\nabla_{0X}\fs+\gamma(X)(\fs),\quad\fs\in
C^\infty(S),\]
we obtain, for the induced connection on $\End(S)$,
\[{\widetilde\nabla}^S_{0X}\ft
={\widetilde\nabla}_{0X}\ft+\gamma(X)\ft-\ft\gamma(X),\] hence
${\widetilde\nabla}^S_0=\nabla$. Although we started from a
Riemannian connection $\nabla_0$ the construction does not
guarantee that $\nabla^S_0$ is also a Riemannian connection. Now
putting
\[\l\fs_1,\fs_2\r'=
X\l\fs_1,\fs_2\r-\l\nabla^S_{0X}\fs_1,\fs_2\r-\l\fs_1,\nabla^S_{0X}\fs_2\r\]
we get a sesquilinear form on $S$, hence
$\l\cdot,\cdot\r'=\omega(X)\l\cdot,\cdot\r$ by (a) of the
Proposition. It is easily seen, that $\omega$ is a (real-valued)
one-form. We can finally put
\[\nabla^S=\nabla^S_0+\frac{1}{2}\omega\] which by (b) of the
Proposition satisfies property (2) and by a simple computation is
seen to be a Riemannian connection with respect to
$\l\cdot,\cdot\r$.
\endproof
\medskip

For $E$ we have found, at least locally, a decomposition $E=S\otimes W$ with 
$\Cl M$ acting irreducibly on $S$. However, there are topological 
obstructions for a global such decomposition to hold. We come back to this 
point later on. However, if $S$ is given globally, $W$ is naturally 
determined by $W=\Hom_{\Cl M_\C}(S,E)$. It is easy to show that this is 
indeed a sub-bundle of the bundle of $\Hom(S,E)$. This gives rise to the 
following definition. 
 
\begin{definition} An oriented Riemannian manifold $M$ of dimension
$m=2n$ is said to be spin$^c$ if there is a complex spinor bundle $S$ over
$M$ with $\Cl M\otimes\C\cong \End(S)$.
\end{definition}

If $M$ is spin$^c$ any spinor bundle $E$ can be written as $E=S\otimes W$ 
with some complex vector bundle $W$. In particular, for any further 
irreducible spinor bundle $S'$ there exists a complex line bundle $L$ with 
$S'= S\otimes L$, viz.\ $L=\Hom_{\Cl M_\C}(S,S')$. Given $S$ we can now make 
it a Dirac bundle by properly choosing a Hermitian structure and a Riemannian 
connection $\nabla$. However, this connection is is only determined up to an 
additional purely imaginary one-form. Most desirable would be a unique 
connection on $S$ induced by the Levi-Civita connection of $M$. Then the 
connection on any further Dirac bundle $S'=S\otimes L$ could be chosen as the 
product connection only depending on the connection on the line bundle $L$. 
To ensure this we need a spin structure for $M$ given by an additional 
structure on $S$.  

We start with the algebraic setting and consider the complex vector space 
$S_0$ of spinors bearing an operation of the real Clifford algebra $\Cl_m$. 
This is not irreducible but depending on the dimension $m=2n=8k+2\ell$ one 
can find an irreducible real subspace of $S_0$. More precisely, there exist 
an antilinear map $\theta_0:S_0\to S_0$ with $\theta_0^2=I_{S_0}$ for 
$\ell=0$ or $3$ and $\theta_0^2=-I_{S_0}$ for $\ell=1$ or $2$, a so-called 
structural map. In the first case $S_0$ carries a real structure, in the 
second case a quaternionic structure. This is obvious if $\ell=0$ or $3$, 
since then $\Cl_m=M_{2^n}(\R)$ is acting irreducibly on $\R^{2^n}$ and the 
complex Clifford algebra and the spinor space $S_0=\C^{2^n}$ are obtained 
therefrom by complexification. Here $\theta$ can be chosen the complex 
conjugation $c:\C^{2^n}\to\C^{2^n}$ taken component-wise. In the other two 
cases we consider the explicit representation of $S_0=\C^{2^n}$, and by 
periodicity may restrict to $k=0$. With $c$ as before and $\tau=i\sigma_2$ we 
now put $\theta_0=\tau\circ c$ if $\ell=1$ and $\theta_0=(\tau\otimes 
\sigma_3)\circ c$, if $\ell=2$. Then $\theta_0$ is the required structural 
map, and in all cases it commutes with the representation of $\Cl_m$. The 
antilinear map $\theta_0:S_0\to S_0$ can also be seen as a linear map 
$\theta_0:S_0\to \bar S_0$, where by $\bar S_0$ we denote the complex vector 
space $S_0$ with scalar multiplication changed to $\lambda\cdot v=\bar\lambda 
v$, $\lambda\in\C$, $v\in S_0$. The representation of $\Cl_m$ on $S_0$ 
induces a representation of $\Cl_m$ on $\bar S_0$, and extending both 
representations to $\Cl_m^\C$ we obtain an element $\theta_0$ of 
$\Hom_{\Cl_m^\C}(S_0,\bar S_0)$ with $\theta_0^2=\pm I_{S_0}$. 

Since $c$ and $\tau$ both depend on a basis of $S_0$ it is in general not 
possible to extend this local construction to a global one on the spinor 
bundle $S$. So at first we will assume a global structural map and afterwards 
will establish sufficient conditions for it existence. 
 
\begin{definition} Let $M$ be an oriented Riemannian
manifold of dimension $m=8n+2\ell$. We say that $M$ carries a spin
structure or that $M$ is spin, if $M$ is spin$^c$ and if the
irreducible complex spinor bundle $S$ allows a structural map
$\theta\in C^\infty\big(\Hom_{\Cl M_\C}(S,\bar S)\big)$ with
$\theta^2=I_S$ or $\theta^2=-I_S$ inducing respectively a real
($\ell=0$ or $3$) or quaternionic ($\ell=1$ oder $2$) structure on
$S$, that is compatible with the complex conjugation of $\Cl
M_\C=\Cl M\otimes\C$.
\end{definition}

{\bf Remark}\quad Equivalently, we may require the existence of a real spinor 
bundle on which the real Clifford bundle $\Cl M$ acts irreducibly on each 
fiber. If $\ell=3$ or $4$ one can choose the fixed-point bundle of $\theta$, 
and conversely the complexified real spinor bundle will define a spin 
structure. 

Of course, any spin manifold is spin$^c$ but the converse does not hold in 
general. We address this question in the next section. Here we only prove the 
following general characterization. 
 
\begin{theorem} Let $M$
be an oriented Riemannian manifold $M$ with spin$^c$ structure given by the
irreducible complex spinor bundle $S$. Then $S$ defines a spin structure,
i.e., allows a global structural map $\theta$ if and only if the vector
bundle $\Hom_{\Cl M_\C}(S,\bar S)$ is trivial.
\end{theorem}

\Proof: We already know that $\Hom_{\Cl M_\C}(S,\bar S)$ is a complex line 
bundle: Each fiber contains a $\Cl(T_pM_\C)$-linear isomorphism 
$\theta_p:S_p\to {\bar S}_p$ and by irreducibility of $S_p$ and ${\bar S}_p$ 
and Schur's Lemma any $\Cl(T_pM_\C)$-linear map $\theta_p':{\bar S}_p\to S_p$ 
satisfies $\theta'_p\circ\theta_p=\lambda I_{S_p}$ for some $\lambda\in\C$. 
In case of a spin structure $\theta$ defines a non-vanishing section in 
$\Hom_{\Cl M_\C}(S,\bar S)$, hence $\Hom_{\Cl M_\C}(S,\bar S)$ is trivial. 
Conversely, if this bundle is trivial and if $\tilde\theta'$ is a 
non-vanishing section, then $\tilde\theta'^2=\lambda I_S$ for some 
non-vanishing map $\lambda\in C^\infty(M,\C)$. But 
\[\lambda(p)\tilde\theta_p(v)
=\tilde\theta_p\circ\tilde\theta_p\circ\tilde\theta_p(v)
=\tilde\theta(\lambda(p)v)=\o{\lambda(p)}\tilde\theta(v),\]
i.e.,
$\lambda\in C^\infty(M)$ is real-valued, and replacing
$\tilde\theta$ by $\theta=|\lambda|^{-1/2}\tilde\theta$ we obtain
a structural map.
\endproof
\medskip

Now given an irreducible complex spinor bundle $S$ and a structural map 
$\theta)$ we can choose a Riemannian structure compatible with Clifford 
multiplication and such that $\theta$ is an isometry. Moreover, we can choose 
a Riemannian connection $\nabla^S$ with properties (1) and (2) uniquely 
determined up to a purely imaginary one-form. If we also require that 
$\nabla^S$ is compatible with $\theta$, i.e.\ 
\[\nabla^S\fs=({\rm id}_{T^\ast
M}\otimes\theta)\nabla^{\bar S}(\theta\circ\fs),\quad \fs\in C^\infty(S),\] 
then such a connection $\nabla^S$ is uniquely determined: 
 
\begin{theorem} If $M$ is a spin manifold of dimension
$m=2n$ with corresponding spinor bundle $S$ and structural map
$\theta$, then:\\ {\rm (a)} On $S$ there exists a Riemannian
structure compatible with Clifford multiplication and with
$\theta$, i.e.\ $\l\theta(\fs_1),\theta(\fs_2)\r=\o{\l
\fs_1,\fs_2\r}$ for $\fs_1,\fs_2\in C^\infty(S)$.\\
{\rm (b)}
There is a unique Riemannian connection $\nabla^S$ with properties
(1) and (2) and compatible with $\theta$.
\end{theorem}

\Proof: (a) We consider $\theta$ as an antilinear map on $S$ and change a 
given Riemannian metric $\l\cdot,\cdot\r'$ with property (1) to 
\[\l\fs_1,\fs_2\r
=\frac{1}{2}\big(\l\fs_1,\fs_2\r'+\o{\l\theta(\fs_1),\theta(\fs_2)\r'}\big).\]
Then the new metric will also be compatible with Clifford
multiplication. Moreover, one has
\[\l\theta(\fs_1),\theta(\fs_2)\r=\frac{1}{2}
\big(\l\theta(\fs_1),\theta(\fs_2)\r'+
\o{\l\theta^2(\fs_1),\theta^2(\fs_2)\r'}\big)
=\o{\l\fs_1,\fs_2\r}.\] In particular,
$\l\theta(\fs_1),\fs_2\r=\pm\l\theta(\fs_2),\fs_1\r$, hence
$\l\theta(\fs_1),\fs_1\r=0$ in the quaternionic case.\\ (b) It
suffices to prove uniqueness. We choose a local orthonormal frame
$\fs_j$ of $S$ (which is a local orthonormal frame of $\bar S$
simultaneously) and the corresponding local connection form
$\omega$. In the real case $S$ is a complexified real spinor
bundle, and we can choose the frame such that
$\theta(\fs_j)=\fs_j$, $j=1,\dots,2^n$. In the quaternionic case
we can choose the frame such that $\fs_{2^{n-1}+j}=\theta(\fs_j)$,
$j=1,\dots,2^{n-1}$. By compatibility of $\nabla^S$ and $\theta$
in the real case we obtain
\[\nabla^{\bar
S}\theta(\fs_j)=\sum_{i=1}^{2^n}\omega_{ji}\fs_i
=\sum_{i=1}^{2^n}\bar\omega_{ji}\fs_i=\theta(\nabla^S\fs_j),\]
i.e.\ $\omega=\bar\omega$. In the quaternionic case we obtain
\begin{eqnarray*}
\nabla^{\bar S}\theta(\fs_j) &=&\sum_{i=1}^{2^n}\omega_{j+2^{n-1},i}\fs_i
=\sum_{i=1}^{2^{n-1}}\bar\omega_{ji}\theta(\fs_i)
+\sum_{i=1}^{2^{n-1}}\bar\omega_{j,i+2^{n-1}}\theta(\fs_{i+2^{n-1}})\\
&=&\sum_{i=1}^{2^{n-1}}\bar\omega_{ji}\fs_{i+2^{n-1}}
-\sum_{i=1}^{2^{n-1}}\bar\omega_{j,i+2^{n-1}}\fs_i =\theta(\nabla^S\fs_j),
\end{eqnarray*}
for $j=1,\dots,2^{n-1}$, i.e.\
\[\omega_{j+2^{n-1},i}=\cases{-\bar\omega_{j,i+2^{n-1}},\quad &
$i=1,\dots,2^{n-1}$,\cr \bar\omega_{j,i-2^{n-1}},\quad &
$i=2^{n-1}+1,\dots,2^n,$\cr}\]
and, in particular,
$\omega_{jj}=\bar\omega_{j+2^{n-1},j+2^{n-1}}$ for
$j=1,\dots,2^{n-1}$.\\ Thus, in both cases addition of a purely
imaginary one-form is prohibited.
\endproof
\medskip

{\bf Examples}\quad 8. Any oriented complex manifold (or, more
generally, an almost-complex manifold) is spin$^c$: Since the
complex cotangent bundle $T^\ast M_\C$) splits orthogonally
$T^\ast
M_\C=(T^cM)^\ast\oplus(\o{T^cM})^\ast=\L^{1,0}M\oplus\L^{0,1}M$,
we can choose $S=\L^\ast \o{T^cM}=\L^{0,\a}M$.

9. However, in general a complex manifold is not spin, e.g.\ it
can be proved that $\C P^n$ is spin if and only if $n$ is odd.

10. Any oriented compact hyper surface $M\subset\R^{2n+1}$ (that
is the boundary of a compact $2n+1$-dimensional submanifold $N$
with boundary) is spin:

Using the matrices $A_j\in M(\C^{2^n})$, $j=1,\dots,2n+1$ the
Clifford multiplication $E_j\cdot v=A_jv$ for $v\in\C^{2^n}=S_0$
and the standard orthonormal frame $E_1,\dots,E_{2n+1}$ of
$\R^{2n+1}$ makes $\R^{2n+1}\times\C^{2^n}$ a (trivial) complex
spinor bundle over $\R^{2n+1}$. If we restrict to $M$ and consider
$TM$ as a subbundle of $TN|_M$ the Clifford modules $\{p\}\times
\C^{2^n}$, $p\in M$, are irreducible $\Cl(T_pM_\C)$ modules, since
we can generate $\Cl_{2n+1}^\C$ by an orthonormal basis
$E'_1(p),\dots,E'_{2n}(p)$ of $T_pM$ and the exterior normal
vector $E'_{2n+1}(p)=X_N(p)$. Therefore, $H=M\times S_0$ defines a
spinor bundle for $M$. Since, moreover, the $E'_j$ as real linear
combinations of the $E_j$ also commute with the structural map
$\theta$ of $S_0$, we even have a spin structure. The grading
operator on $H$ is defined by $\epsilon=-iX_N\cdot$ with the
exterior normal vector field $X_N$ at $M$. If $M=S^{2n}$, we have
$X_N(x)=\sum_{k=1}^{2n+1}x_kE_k$ and $H^{n (\mod 2)}$, the bundles
of half-spinors are non-trivial smooth vector bundles.

Special examples are oriented compact surfaces $T_g$ in $\R^3$ or spheres 
$S^{2n}$ in $\R^{2n+1}$. The former allow $2^{2g}$ different spin structures 
whereas there is only one spin structure on $S^{2n}$. To see this we need the 
following result. 
 
\begin{theorem} Let $M$ be a connected oriented
Riemannian manifold. If $M$ carries a spin$^c$ structure, then all
of the non-equivalent spin$^c$ structures are parametrized by
$H^2(M,\Z)$. If moreover $M$ is spin, then all of the different
spin structures are parametrized by $H^1(M,\Z_2)\cong
\Hom\big(\pi_1(M),\Z_2\big)$. In particular, $M$ allows at most
one spin structure if $M$ is simply connected.
\end{theorem}

\Proof: Starting from a irreducible complex spinor bundle $S$, any further 
irreducible complex spinor bundle on $M$ is of the form $S'=S\otimes L$ where 
$L=\Hom_{\Cl M_\C}(S,S')$. If $S'$ and $S''=S(E)\otimes L'$ are isomorphic as 
spinor bundles, i.e.\ determine equivalent spin$^c$ structures, there is a 
$\Phi\in\Iso_{\Cl M_\C}(S',S'')$, and so $L\cong L'$. This shows that 
$H^2(M,\Z)$ acts transitively on the set of different spin$^c$ structures. 
Now for $\Hom_{\Cl M_\C}(S',\bar S')$ we obtain 
\begin{eqnarray*}
\Hom_{\Cl M_\C}(S',\bar S')&\cong& \Hom_{\Cl M_\C}(S\otimes_\C L,\bar
S\otimes_\C \bar L)\cong \Hom\big(L,\Hom_{\Cl M_\C}(S,\bar S\otimes_\C \bar
L)\big)\\ &\cong& L^\ast\otimes_\C \Hom_{\Cl M_\C}(S,\bar S \otimes_\C \bar
L)\cong L^\ast\otimes_\C \Hom_{\Cl M_\C}(S,\bar S \otimes_\C \bar L\cr
&\cong& L^\ast\otimes_\C \bar L\cong \Hom(L,\bar L),
\end{eqnarray*}
if $\Hom_{\Cl M_\C}(S,\bar S)$ is trivial. Therefore, there is a
structural map on $S'$ if and only if $L^\ast\otimes \bar L\cong
{\bar L}^2$ is trivial. If $H^2(M,\Z)$ has no 2-torsion, $\bar L$
has to be trivial, too, and likewise $L$. In any case different
spin structures are classified by isomorphy classes of real line
bundles, i.e., by $H^1(M,\Z_2)$; cf.\ [Kar].
\endproof
\medskip

{\bf Remarks} 1. We always started with the Clifford bundle of the tangent 
bundle. Only with literate changes we can start with a real Riemannian vector 
bundle $E$ of even rank. A spin$^c$ structure is then given by a complex 
spinor bundle $S(E)$ with $\Cl^\C(E)$ acting irreducibly on the fibers, and a 
spin structure by an additional structural map compatible with Clifford 
multiplication. If $E$ comes with a Riemannian connection $\nabla^E$ there is 
unique connection $\nabla^{\Cl(E)}$ on $\Cl(E)$ and in the spin case a unique 
Riemannian connection $\nabla^{S(E)}$ on $S(E)$ that satisfy properties (1) 
and (2) and 
\[\nabla^{S(E)}(\sigma\cdot\fs) = \nabla^{\Cl(E)}(\sigma)\cdot \fs +
\sigma \cdot \big(\nabla^{S(E)}\fs\big),\]
for $\sigma \in
C^\infty(\Cl(E))$, $\fs \in C^\infty(S(E))$.

2. On an oriented Riemannian vector bundle $E$ of odd rank $m=2n+1$ (in
particular, on an odd-dimensional Riemannian manifold) spin$^c$ or spin
structures can be defined, too. Here a spin$^c$ structure is given by a
complex spinor bundle $S(E)$, on which $\Cl^\C(E)$ acts irreducibly, and
where for each oriented orthonormal frame $e_1(p),\dots,e_m(p)$ of $E_p$ the
element $i^{n+1}e_1(p)\cdots e_m(p)$ acts as $I_{E_p}$.
\bigskip
\bigskip
\section{\bf Spin groups and principal bundles}
\bigskip
There are topological obstructions for a spin$^c$ or a spin
structure to exist on a manifold $M$. We know that if $M$ is
spin$^c$ and $S$ an irreducible complex spinor bundle structure
then $M$ is spin if and only if $L=\Hom_{\Cl M_\C}(S,\bar S)$ is
trivial. Now if $M$ is simply connected this can be decided by
computing a topological invariant. It is well known (cf.\ [Sdr2])
that $L$ is trivial if and only if the first Chern class $c_1(L)$
vanishes. But this does not apply in general if $M$ is not simply
connected. Then the obstructions are better expressed in terms of
the so-called second Stiefel-Whitney class $w_2(TM)$, an element
of $H^2(M,\Z_2)$ (cf.\ [Hae]). This is a cohomology class with
coefficients in $\Z_2=\{\pm 1\}$, and can be represented by lifts
of cocycles of $SO(n)$-valued transition maps to the covering
group $Spin(n)$.
\medskip

At this point we have to digress and take a closer look at the
covering group $Spin(n)$ of $SO(n)$. Here again Clifford algebras
are the appropriate tool to generalize classical constructions. We
first inspect how Clifford algebras help represent orthogonal
transformations. It is well known that $S^3\subset\H$ is the
two-fold simply connected covering of the Lie group $SO(3)$.
Identifying $\R^3$ with ${\rm Im}~\H=\{is+jt+ku\in\H\mid
s,t,u\in\R\}$ an element $x\in S^3=\{y\in\H\mid |y|^2=\bar y
y=1\}$ acts on $\R^3$ by
\[\Ad_x(v) = xvx^{-1}=xv\bar x,\quad v\in\R^3.\]
Note that $x$ and $-x$ define the same element of
$SO(3)$. More generally one could use any
$x\in\H^\ast=\H\setminus\{0\}$ since $\Ad_x=\Ad_{x/|x|}$.

To find the covering group of $SO(n)$ for $n\ge 4$ or of $SO(E)$
for a Euclidean vector space $E$ we start from the regular group
$G\Cl(E)$ of invertible elements of the algebra $\Cl(E)$. For
$x\in E\setminus \{0\}\subset G\Cl(E)$ and $v\in E\subset\Cl(E)$
we have $x\cdot v+v\cdot x=-2\l x,v\r 1$, hence
\[-\Ad_x(v)=v-2\frac{\l x,v\r}{\l x,x\r} x.\]
From a geometric point of view this is the reflection at the
hyperplane perpendicular to $x$. Using the involution $\alpha$
(that induces the grading $\Cl(E)=\Cl(E)^0\oplus\Cl(E)^1$) we pass
over to the ``twisted'' adjoint representation on $E$ given by
\[\widetilde{\Ad}_x(v)=\alpha(x)vx^{-1}\]
which is is naturally defined on the Clifford group
\[\Gamma(E)=\{x\in G\Cl(E)\mid \alpha(x)vx^{-1} \in E~\hbox{ for all
}v \in E\}.\]  

\begin{propo} The twisted adjoint representation $\widetilde{Ad} : \Gamma(E) \to \Aut(E)$
is a homomorphism of groups and induces an exact sequence
\[1 \to \R^\ast \to \Gamma(E) ~{\buildrel \scriptstyle{\widetilde{\Ad}}
\over \longrightarrow} ~O(E) \to 1.\]
Any $x \in\Gamma(E)$ can be
written as $x = v_1 \cdots v_k$, $v_i \in E$, $v_i \ne 0$, $i =
1,\dots,k$.
\end{propo}

\Proof: Obviously, $\widetilde{\Ad}$ is a homomorphism. Next we show that 
$x\in\R^\ast =\R\setminus\{0\}\subset \Gamma(E)$ if $\alpha(x)v=vx$ for all 
$v\in E$ or equivalently if this holds elements $v$ of an orthonormal basis 
$(e_i)_{1\le i \le n}$ of $E$. To this end we write $x = x^0 + x^1 \in 
\Cl(E)^0 \oplus \Cl(E)^1$ with $x^0 = a_i^0 + e_ib_i^1$ and 
$x^1=a^1_i+e_ib_i^0$, where $a_i^j$ and $b_i^j$ are of degree $j$ (mod 2) and 
both do not contain $e_i$. Then we get 
\[\alpha(x)e_i=(x^0-x^1)e_i=e_i(a^0_i+a^1_i)+b^1_i+b^0_i\]
and
\[e_ix=e_i(x^0+x^1)=e_i(a^0_i+a^1_i)-b^1_i-b^0_i,\]
which entails $b^0_i=b^1_i=0$, i.e.\ $x\in\R^\ast$.

Since $O(E)$ is generated by reflections it is at least contained
in the image of $\widetilde\Ad$.

It remains to show $\widetilde\Ad\big(\Gamma(E)\big)\subset O(E)$,
i.e.\ $|\widetilde\Ad_x(v)|=|v|$ for $v\in E$. To prove this we
consider the anti-automorphism of $\Cl(E)$ induced by
\[x = v_1 \cdots v_k \mapsto x^t = v_k \cdots v_1\]
and the anti-automorphism
\[\Cl(E)\ni x \mapsto \b x =
\alpha(x^t)=\big(\alpha(x)\big)^t\in\Cl(E)\] which allows to
extend the quadratic form $E\ni v\mapsto v \cdot \b v = \l v,v\r 1
= |v|^2 1\in\Cl(E)$ to the so-called spinor norm
\[\Cl(E) \ni x\mapsto N(x) = x \cdot \b x\in \Cl(E)\]
of $\Cl(E)$. Since the
anti-automorphisms leave $\Gamma(E)$ invariant, we have
$N\big(\Gamma(E)\big) \subset \Gamma(E)$. Actually
$N\big(\Gamma(E)\big)\subset\R^\ast$, because
\begin{eqnarray*}
\widetilde\Ad_{N(\b x)}(v) &=&
\alpha\big(\alpha(x^t)x\big)v\big(\alpha(x^t)x\big)^{-1}
=x^t\alpha(x)vx^{-1}\alpha(x^{-1})^t\\
&=&\big(\alpha(x^{-1})\alpha(x)vx^{-1}x\big)^t=v.
\end{eqnarray*}
Now $N|_{\Gamma(E)}$ is a homomorphism of groups, since
\[N(xy) = xy\o{xy} = xy\alpha(y^t)\alpha(x^t) = xN(y)\alpha(x^t) =
N(x)N(y),\] as $N(\Gamma(E))\subset\R^\ast$. In particular,
\[N\big(\alpha(x)vx^{-1}\big)= N\big(\alpha(x)\big)N(v)N(x)^{-1} =
N(v)N\big(\alpha(x)\big)N(x)^{-1} = N(v),\] since
$N\big(\alpha(x)\big) = \alpha(x)x^t =\alpha\big(N(x)\big)=N(x)
\in \R^\ast$, and we conclude
\[|\widetilde{\Ad}_x(v)|^2 = |\alpha(x)vx^{-1}|^2 = |v|^2,\]
i.e.\ $\widetilde{\Ad}_x \in O(E)$.
\endproof

\begin{definition} We put $Pin(E)=N^{-1}(1)\cap \Gamma(E)$
and define the spin group of the Euclidean vector space $E$ by
$Spin(E) = Pin(E) \cap \Cl(E)^0$. In the case $E = \R^n$ with its
standard inner product we write $Spin(n)$ instead of $Spin(\R^n)$.
\end{definition}

{\bf{Remarks}}\quad 1. The group $Spin(E)$ is compact, in fact a Lie group as
a closed subgroup of the group of invertibles of the algebra $\Cl^0(E)$.

2. Of course, $Pin(E)$ and $Spin(E)$ both depend on the Euclidean structure.
More generally, one can also define $Spin(E,Q)$ for a real vector space $E$
and a non-degenerate quadratic form $Q$.

3. One has ${\rm Pin}(E) = \{v_1 \cdots v_{k} \in \Cl(E) \mid v_i \in E, \l
v_i,v_i\r  = 1,~i=1,\dots,k \}$ and $Spin(E) = \{v_1 \cdots v_{2k} \in \Cl(E)
\mid v_i \in E, \l v_i,v_i\r  = 1,~i=1,\dots,2k \}$.
\medskip

\noindent {\bf Corollary}~ {\it The groups $Pin(E)$ and $Spin(E)$ fit into 
the following exact sequences 
\[1 \to \Z_2 \to Pin(E) \to O(E) \to 1\]
\[1 \to \Z_2 \to Spin(E)  \to SO(E) \to 1.\]
In particular,
\[1 \to \Z_2 \to Spin(n) \to SO(n) \to 1\]
is exact, i.e., $Spin(n)$ is a non-trivial two-sheeted covering of
$SO(n)$. For $n\ge 3$ it is simply connected, i.e.\ the universal
covering group of $SO(n)$.}
\smallskip

\Proof: Given $x \in \Gamma(E)$ and $\lambda = 1/\sqrt{N(x)}$ one
has $\lambda x \in Pin(E)$ hence
\[\widetilde{\Ad}|_{Pin(E)} : Pin(E) \to O(E)\]
is onto and
\[\Ker\widetilde{\Ad}|_{Pin(E)} = \{ \lambda \in \R^\ast
\mid N(\lambda) = \lambda^2 = 1 \} \cong \Z_2.\] Any element of
$SO(E)$ may be written as $\widetilde{\Ad}_{v_1} \cdots
\widetilde{\Ad}_{v_{2k}}$ hence
\[\rho =
\widetilde{\Ad}|_{\Gamma(E) \cap\Cl(E)^0} : \Gamma(E) \cap
\Cl(E)^0 \to SO(E) \]
is onto with $\Ker \rho = \R^\ast$. Now the
restriction to $Spin(E)$ yields the analogous exact sequence. To
prove the last assertion we only have to find a continuous path
connecting $+1$ and $-1$ in $Spin(n)$. To this end we choose
$e_1,e_2 \in \R^n$ with $e_1 \bot e_2,\, |e_i| = 1$, and
\begin{eqnarray*}
c(t) & = & \exp(2\pi te_1\cdot e_2) =\cos 2\pi t + e_1 \cdot e_2 \sin 2\pi t
\\ & = & (e_1 \cos \pi t + e_2 \sin \pi t)\cdot (-e_1 \cos \pi t + e_2 \sin \pi
t),
\end{eqnarray*}
for $0 \le t \le \frac{1}{2}$. Thus the covering is non-trivial
and $Spin(n)$ is connected. For $n \ge 3$ it is also simply
connected by the classic topological result $\pi_1\big(SO(n)\big)
= \Z_2$, $n \ge 3$.
\endproof
\medskip

If $E_\C$ is the complexification of $E$ with $\C$-linear extension $Q_\C$ of
$Q$, then $\Cl(E_\C,Q_\C)$ and $\Cl(E,Q)\otimes\C$ are isomorphic. We put
$\alpha(x\otimes z)=\alpha(x)\otimes z$ and $(x\otimes z)^t=
x^t\otimes\bar{z}$ and with $\bar{}~$ and $N$ as before we also define
$Pin^c(E)$ and the group $Spin^c(E)\subset \Cl^0(E,Q)\otimes\C$. The latter
is isomorphic with $Spin(E)\times S^1/\Z_2$ where $\Z_2=\{(1,1),(-1,-1)\}$.
If $E=\R^n$ we simply denote it by $Spin^c(n)$.
\medskip
The group $Spin^c(E)$ is also compact and fits into the exact
sequences
\[1 \to S^1~ \to~ Spin^c(E)~ {\buildrel \rho_0\over \longrightarrow}~ SO(E) \to 1\]
\[1 \to Spin(E)~ \to~ Spin^c(E)~
{\buildrel \rho_1\over\longrightarrow}~ S^1 \to 1,\]
where the
left hand homomorphisms are canonical inclusions and the right
hand ones are defined by $\rho_0([(x,z)])=\rho(x)$ and
$\rho_1([(x,z)])=z^2$, $(x,z)\in Spin(E)\times S^1$, respectively.
\medskip

Usually, spin and spin$^c$ structures are defined with the help of
corresponding principal bundles; cf.\ [BH] and [Mil]. One starts
with the orthonormal frame bundle $P_{SO(m)}$ of the tangent
bundle of an $m$-dimensional oriented Riemannian manifold $M$ (or
of an oriented Riemannian vector bundle of rank $m$). A spin
structure for $M$ consists of a principal bundle $P_{Spin(m)}$
with structure group $Spin(m)$ and a two-sheeted covering
\[\xi :
P_{Spin(m)} \to P_{SO(m)}\quad\hbox{ with }\quad\xi(p g) =
\xi(p)\rho_0(g),~ p\in P_{Spin(m)},~ g \in Spin(m),\] where
$\rho_0 : Spin(m) \to SO(m)$ is the standard covering. A spin$^c$
structure is given by a principal bundle $P_{Spin^c(m)}$ and a map
\[\xi : P_{Spin^c(m)} \to P_{SO(m)}\quad\hbox{ with }\quad \xi(p g)
= \xi(p)\rho_0(g),~ p\in P_{Spin^c(m)},~ g \in Spin^c(m)\]
where
$\rho_0: Spin^c(m)\to SO(m)$ is again the standard map.

To show that this approach is equivalent with the one presented so
far one has to go two ways. A spinor bundle can be obtained as an
associated bundle: If $F$ is a real or a complex vector space,
which is also a $\Cl_m$-module or a $\Cl_m^\C$-module with
compatible inner product, representations $\rho : Spin(m) \to
SO(F)$ or $\rho : Spin^c(m)\to U(F)$ will be induced by
left-multiplication with elements of $Spin(m) \subset \Cl^0_m$ or
$Spin^c(m)\subset \Cl_m^0\otimes\C$, respectively. Then $S =
P_{Spin(n)} \times_\rho F$ is a real or a complex spinor bundle,
which moreover is irreducible if $F=S_0$ the space of spinors.
\\If on the other hand a spin structure is given by an irreducible
complex spinor bundle $S$ the corresponding principal bundles can
be recovered as follows. First recall that $P_{SO(m)}$ can be
considered as the subset of $\Hom(M\times\R^m,TM)$ that consists
of all orientation preserving isometries $f_p:\R^m\to T_pM$, $p\in
M$. Then we define $P_{Spin(m)}$ and $P_{Spin^c(m)}$ to be
appropriate subsets of $\Hom(M\times S_0,S)$. In the second case
it consists of all isometries $\phi_p:S_0\to S_p$ that respect the
decompositions $S_0^0\oplus S_0^1$ and $S_p^0\oplus S_p^1$ and
satisfy $\phi_p(v\cdot\phi_p^{-1})\in T_pM\subset\Cl(T_pM_\C)$ for
all $v\in\R^m\subset\Cl_m^\C=\End(S_0)$. In the first case we
additionally require that these isometries respect the real or
quaternionic structure. The map $\xi: P_{Spin^c(m)}\to P_{SO(m)}$
is now defined by $\xi(\Phi_p)=\Ad(\Phi_p)$. Then one has
\[\xi(\Phi_pg)=\Ad(\Phi_pg)=\Ad(\Phi_p)\circ\Ad(g)=\xi(\Phi_p)\rho_0(g).\]
This action from the right is transitive, since for
$\Phi_p,\Phi'_p\in P_{Spin^c(m)}$ one has
$\Phi_p=\Phi'_p\circ((\Phi'_p)^{-1}\circ\Phi_p)$ and by definition
$x=(\Phi'_p)^{-1}\circ\Phi_p \in\Cl_m^0\otimes\C$ as well as
$N(x)=1$, hence $x\in Spin^c(m)$. Here $N(x)=1$ does hold, since
$x$ is unitary and since $\bar x=x^\ast$ for
$x\in\Cl_m\otimes\C=\End(S_0)$ as $\alpha(v)=-v=v^\ast$ for $v\in
\R^m$. In the spin case $\xi$ is defined likewise and obviously
such an element $x$ belongs to $Spin(m)$.

Finally, we can come back to the topological obstructions that
decide upon spin$^c$ or spin structures. A spin$^c$ structure can be
supplied if and only if $w_2(TM)$ is the $\mod 2$-reduction of
some integral cohomology class (or, what amounts to the same, if
the integral Stiefel-Whitney class $W_3(TM)$ vanishes). A spin
structure exists if and only if $w_2(TM)=1$. We refer to [Kar2],
where the first assertion is proved explicitly and the second one
implicitly -- in the case of a spin structure in the notation of
[Kar2] one has to replace $\C^\ast$ by $\R^\ast$, which makes
$l_{1\ast}$ automatically an isomorphism. These conditions can be
checked combinatorically. We refer to [Gil] and [LM] for some
specific computations. In particular, it is proved that $w_2(T\C
P^n)$ and $w_2(T\R P^{2n+1})$ only vanish if $n$ is odd. Using
deeper results of algebraic topology one can show that any compact
oriented 3-manifold is spin (since according to E.\ Stiefel it is
parallizable) and that any compact oriented 4-manifold is spin$^c$
(according to a theorem of Whitney; cf.\ [HH]).
\bigskip
\bigskip
\section{\bf The geometric Dirac operators}
\bigskip
Now we want to look more closely at some Dirac operators. First we
consider the special case $M=\R^m$ with its standard metric and
the global orthonormal frame $E_j = \frac{\d}{\d x_j}$,
$j=1,\dots,m$, of $T\R^m$. If $V$ is an $n$-dimensional
$\Cl_m$-module defined by an algebra homomorphism $\rho : \Cl_m
\to \End(V)$, say $\rho(e_j) = A_j$, $\rho(1) = I_V$, and
$(v_i)_{1\le i \le n}$ is a basis of $V$ a global frame on
$E=\R^m\times V$ is given by $\fs_i(p) = (p,v_i)$, $p\in \R^m$,
$i=1,\dots,n$. Let $\nabla$ denote a flat connection on $E$, i.e.\
$\omega\equiv 0$ with respect to the frame $\fs_i$, hence $\nabla
f \fs_i = d f \otimes \fs_i$ for $f\in C^\infty(\R^m)$. Then
Clifford multiplication $E_j \cdot \fs_i$ is given by
\[(E_j\cdot \fs_i) (p) = \big(p,A_j(v_i)\big)\,,\] hence
\[D\Big(\sum_{i=1}^n f_i \fs_i\Big)(p) = \sum_{j=1}^m \sum_{i=1}^n
E_j \cdot \nabla_{E_j} (f_i\fs_i)(p) = \Big(p,\sum_{j=1}^m
\sum_{i=1}^n \frac{\d}{\d x_j} f_i(p) A_j(v_i)\Big)\,.\]

Let $A_j$ also denote the matrix with respect to the basis
$(v_i)$. Then a (local) representation of $D$ is given by
\[D f = \sum_{j=1}^n A_j \frac{\partial}{\partial x_j} f\,,\]
where $f=(f_1,\dots,f_n)^T$. In particular
\[D^2 f =
\sum_{j,k=1}^m A_j A_k \frac{\d}{\d x_j} \frac{\d}{\d x_k} f =
\Delta\otimes I_n f = - \sum_{j=1}^m \frac{\d^2}{\d x_j^2}f\,,\]
since
\[A_jA_k + A_k A_j = \rho (e_j\cdot e_k + e_k \cdot e_j) =
\rho (-2\delta_{jk}) = - 2 \delta_{ik} I_V.\]
Thus $D$ is a
square-root of $\Delta$. In cases $m=1,2$ we have the following
classical operators.
\medskip

{\bf{Examples}}\quad 11. If $m=1$, i.e. $\Cl_1=\C$, we choose
$V=\C\cong \R^2$ with $\rho(e_1) = i$ and obtain $D= i\frac{\d}{\d
x}$.

12. If $m=2$, i.e.\ $\Cl_2 = \H$, we choose $V=\H$ with $\rho(e_1)
= i$, $\rho(e_2) = j$ and get a grading $\Cl_2 = \Cl_2^0 \oplus
\Cl_2^1 \cong \C \oplus \C$ by
\begin{eqnarray*}
&&\Cl_2^0 \ni u + v e_2 e_1 \mapsto u+ iv \in\C\\
&&\Cl_2^1 \ni u e_1 + v e_2
\mapsto u+iv \in\C.
\end{eqnarray*}
Identifying $E=\R^2 \times V$ and $\C \times (\C \oplus \C)$ the
Dirac operator $D=i \frac{\d}{\d x_1} + j \frac{\d}{\d x_2}$
becomes
\[\frac{1}{2}D(f\oplus g) = - \frac{\d}{\d z} g \oplus {\d
\over \d \bar z} f\,,\quad f,g \in C^\infty (\C)\,.\]
If we write
$D=D^0\oplus D^1$ with $D^j : C^\infty(\R^2 \times \Cl_2^j)\to
C^\infty\big(\R^2 \times \Cl_2^{j+1~(\mod 2)}\big)$ then ${1\over
2}D^0$ is just the Cauchy-Riemann operator $\bar \d = {\d \over \d
\bar z}$ which is studied in the theory of complex functions.
\medskip

If $M$ is a (compact) oriented Riemannian manifold there are several Dirac
operators related to additional geometric structures. We cannot go into the
analytic properties of these Dirac operators; cf.\ [Gil], [LM], or [Sdr1,2].
We only note that they are symmetric (elliptic) differential operators.

If $M$ is of even dimension $m=2k$ we have a global section
$\omega\in C^\infty(\Cl M_\C)$ which is locally given by
\[\omega = i^k E_1 \cdot\cdots\cdot E_m\]
with respect to an oriented orthonormal frame $(E_i)_{1\le i\le
m}$ of $TM$. Obviously, one has $\omega^2 = 1$ and $\omega \cdot X
= - X\cdot \omega$ for $X\in C^\infty(TM)$. Since $\omega$ does
not depend on the local frame we may assume
$\o{\nabla}\nolimits_{E_i} E_j(p) = 0$ at a fixed point $p\in M$
and conclude that
\[\o{\nabla}\nolimits_{E_i} \omega(p) = i^\ell
\sum_{j=1}^m E_1 \cdots \nabla_{E_i} E_j \cdots E_m(p) = 0,\]
hence $\o{\nabla} \omega = 0$. Using $\omega$ any Dirac bundle $E$
on $M$ will be graded by $E^0 = (1+\omega)\cdot E$ and $E^1 =
(1-\omega)\cdot E$. For $\fs\in C^\infty(E^j)$ and $X\in
C^\infty(TM)$ we then obtain
\[\nabla_X \fs = (-1)^j\nabla_X
(\omega\cdot \fs) = (-1)^j\big(({\o{\nabla}}_X \omega) \cdot \fs +
\omega\cdot \nabla_X\fs\big) =(-1)^j \omega \cdot \nabla_X \fs,\]
i.e.\ $\nabla_X\fs ={1\over
2}\big(1+(-1)^j\omega\big)\nabla_X\fs\in C^\infty(E^j)$, and
\[X\cdot \fs = (-1)^jX\cdot\omega\cdot \fs= (-1)^{j+1}\omega \cdot
X\cdot \fs,\] hence $X\cdot \fs ={1\over 2}\big(1+(-1)^{j+1}\omega\big)\cdot 
X\cdot\fs\in C^\infty(E^{j+1 \mod 2})$. Since 
$\l\omega\cdot\fs,\ft\r=\l\fs,\omega\cdot\ft\r$ for $\fs,\ft\in C^\infty(E)$, 
the decomposition $E=E^0\oplus E^1$ is orthogonal. This gives rise to the 
following definition:  
\medskip

\begin{definition} Let $E$ be a Dirac bundle on $M$. An
orthogonal decomposition $E=E^0\oplus E^1$ is called admissible if

{\rm (1)} $\Cl^i(M) E^j \subset E^{i+j \mod 2}$,

{\rm (2)} $\nabla_X \fs\in C^\infty(E^i)$ for $\fs\in C^\infty(E^i)$, $X\in
C^\infty(TM)$.
\end{definition}

Now given an admissible Dirac bundle $E=E^0\oplus E^1$ the corresponding 
Dirac operator induces first order differential operators 
\[D^i : C_0^\infty(E^i) \to C_0^\infty(E^{i+1 \mod 2}),\quad i=0,1.\]
If $M$ is compact $D$ and $D^j$ are elliptic and extend to bounded
operators on appropriate Sobolev space. These extensions are
Fredholm operators, i.e.\ have an index
\[\ind D=\dim\ker D-\dim\coker D\in\Z.\]
Of course, $\ind D=0$ but
\[\ind D^0=\dim\ker D^0-\dim\ker D^1\]
turns out to be an interesting geometric invariant. We already met the 
K\"{a}hler-Dirac operator $D=d+\delta$. Since the Dirac bundle $(\L^\ast 
M=\L^{\rm ev} M\oplus\L^{\rm odd} M,\o{\nabla})$ is admissible, we obtain 
$D=D^0\oplus D^1$ and $\ind D^0=\chi(M)$, the Euler characteristic of $M$. 
There is a different admissible decomposition of $\L^\ast M$ given by 
$\omega$. If $M$ is of dimension $m=4k$ the index of the corresponding Dirac 
operator is just the signature of $M$; cf.\ [Gil], [LM], or [Sdr2]. We can 
now, finally, define the Spin-Dirac or Atiyah-Singer operator. 
 
\begin{definition} Let $M$ be a compact spin manifold of dimension
$m=2n$ with spinor bundle $S$, and $\nabla^S$ the unique connection on $S$
that is induced by the Levi-Civita connection. The Dirac operator associated
with the Dirac bundle $S$ is called the Spin-Dirac operator or
Atiyah-Singer operator and will be denoted by $D_{AS}$.
\end{definition}

The Dirac bundle $S = S^0 \oplus S^1$ with decomposition induced by $\omega$ 
is admissible, i.e. $D_{AS}=D^0_{AS}\oplus D^1_{AS}$. The operator 
$\nd=D_{AS}^0$ is also often called the Spin-Dirac operator. 

Its index $\ind {\nd}=\hat A(M)$ is called the $\hat A$-genus of
$M$. It has topological significance which is expressed by the
famous Atiyah-Singer index theorem:
\[\ind {\nd}=\int_m\hat A(TM)\]
where $\hat A(TM)\in H^m(M,\R)$ is the cohomology class first indroduced by 
F.\ Hirzebruch in 1954; cf.\ [Gil] for a detailed history of the subject 
matter. 

We finally take a closer look on the Dirac-Laplace operator $D^2$
and on its relation to the curvature tensor $\curv (\nabla)$ which
is defined by
\[\curv(\nabla)(X,Y)=[\nabla_X,\nabla_Y]-\nabla_{[X,Y]}\]
for vector fields $X$ and $Y$. Here the brackets denote the
respective commutators. Note that the curvature tensor
$\curv(\nabla)(X,Y)$ of a Riemannian connection is skew-adjoint,
\[\l\curv(X,Y)\fs,\ft\r=-\l\fs,\curv(X,Y)\ft\r\hfill(R)\]
and that
\[\curv(\nabla)(X,Y)(\sigma\cdot\fs)=\curv(\o{\nabla})(X,Y)(\sigma)\cdot\fs
+\sigma\cdot\curv(\nabla)(X,Y)\fs,\hfill(D)\] if $\nabla$
satisfies property (2) of a Dirac triple.

Recall that the second covariant derivative
\[\nabla_{X,Y}^2 : C^\infty(E) \to C^\infty(E)\,,\]
is defined for $X,Y\in C^\infty(TM)$ by
\[\nabla_{X,Y}^2 \fs =
\nabla_X \nabla_Y \fs - \nabla_{\o{\nabla}\nolimits_X Y} \fs \,,
\quad \fs\in C^\infty(E)\,,\] and the curvature tensor of $\nabla$
is given by
\[\curv (\nabla)(X,Y) = \nabla^2_{X,Y} -
\nabla_{Y,X}^2.\] Since $\nabla_{X,Y}^2\fs(p)$ only depends on
$X_p$ and $Y_p$ it makes $\nabla_{\cdot,\cdot}^2$ and $\curv
(\nabla)$ tensors with values in $E_p$. The Bochner-Laplace
operator of the connection $\nabla$ is defined as
\[\nabla^\ast
\nabla \fs = - \tr (\nabla_{\cdot,\cdot}^2 \fs)\,,\quad \fs\in
C^\infty(E)\,,\] i.e.\ as
\[\nabla^\ast\nabla \fs = - \sum_{j=1}^m
\nabla_{E_j,E_j}^2 \fs \]
when computed in some local orthonormal
frame $(E_j)_{1\le j\le m}$. This definition does not depend on
the chosen frame. Moreover, $\nabla^\ast\nabla : C^\infty(E) \to
C^\infty(E)$ is a second order (elliptic) differential operator.

Using the smooth section ${{\cal R}}\in C^\infty(\Hom (E,E))$
given by
\[{{\cal R}}(\fs) = {1\over 2} \sum_{j,k=1}^m E_j \cdot
E_k \cdot \curv (\nabla) (E_j,E_k)(\fs)\,\] we obtain the following 
fundamental result.  

\begin{theorem}{\rm (Bochner-Weitzenb\"{o}ck)}\quad Let
$E$ be a Dirac bundle over $M$ with associated Dirac operator $D$.
Then the Dirac-Laplace operator satisfies $D^2$
\[D^2 = \nabla^\ast \nabla + {\cal R}\,.\]
\end{theorem}

\Proof: With the frame $(E_j)_{1\le j\le m}$ at $p$ as above we have
\begin{eqnarray*} D^2 &=&\sum_{j,k=1}^m E_j \cdot \nabla_{E_j} (E_k \cdot
\nabla_{E_k}) = \sum_{j,k=1}^m E_j\cdot E_k \cdot \nabla_{E_j}
\nabla_{E_k}\\ &=& \sum_{j,k=1}^m E_j \cdot E_k \cdot
\nabla_{E_j,E_k}^2 \\ &=& - \sum_{j=1}^m \nabla_{E_j,E_j}^2 +
\sum_{1\le j<k\le m} E_j \cdot E_k\cdot (\nabla_{E_j,E_k}^2 -
\nabla_{E_k,E_j}^2)\\ &=&\nabla^\ast \nabla +{\cal R}\,.
\end{eqnarray*}
\smallskip

A simple but important application is the following vanishing theorem.
\medskip

\noindent
{\bf Corollary}~ {\it If $M$ is compact and connected and ${\cal
R}(p)$ positive semi-definite for all $p\in M$ and positive
definite for at least one $p$, then the differential equation
$D^2\fs =0$ has only the trivial solution $\fs = 0$, i.e.\ there
are no non-trivial harmonic sections.}
\smallskip

\Proof: For a fixed point $p\in M$ one can choose $(E_j)_{1\le
j\le m}$ such that ${\o{\nabla}}_{E_i}E_j(p) = 0$, and given
sections $\fs,\ft\in C^\infty(E)$ there is a vector field $X$ with
\[\l X,Y\r  = \l  \nabla_Y \fs,\ft\r_E\,,\quad Y\in
C^\infty(TM)\,.\]
These data help to prove that
\begin{eqnarray*}
\l \nabla^\ast\nabla \fs,\ft\r (p) &=& - \sum_{j=1}^m \l \nabla_{E_j}
\nabla_{E_j} \fs, \ft\r (p)\\ &=& - \sum_{j=1}^m \big( E_j \l \nabla_{E_j}
\fs,\ft \r  - \l \nabla_{E_j} \fs,\nabla_{E_j} \ft\r \big)(p)\\ &=& - {\rm
div}\; X(p) + \l \nabla \fs,\nabla \ft\r(p)\,.
\end{eqnarray*}
When integrated over $M$, by Gau{\ss}' theorem, the divergence term
does not occur and we obtain
\[0\le\int_M\l {\cal R}(\fs),\fs\r=-\int_M\l
\nabla^\ast\nabla\fs,\fs\r=-\int_M\l\nabla\fs,\nabla\fs\r \le 0,\]
if $D^2\fs=0$. Therefore,  $\nabla\fs=0$ and $\|\fs\|$ is
constant, since $\nabla$ is Riemannian. Assuming $s(p)\ne 0$ and
${\cal R}(p)$ positive definite gives $\int_M\l{\cal
R}(\fs),\fs\r_E>0$, which cannot hold.
\endproof
\medskip 

There are a lot of special cases of the Bochner-Weitzenb\"{o}ck formula. The 
Bochner-Weitzenb\"{o}ck formula for the Laplace operator can already be found in 
Weitzenb\"{o}ck's monograph ``Invariantentheorie'' of 1923. It has been 
rediscovered and applied in 1946 by S.\ Bochner [Boc]. Here we only consider 
one special case and deduce a special vanishing result. 

\begin{theorem} If $M$ is a spin manifold with spinor
bundle $S$ and connection $\nabla^S$, then the Spin-Dirac-Laplace
operator $D_{AS}^2$ and the Bochner-Laplace operator
$\nabla^{S\a}\nabla^S$ are related by
\[D_{AS}^2 = \nabla^{S\a}\nabla^S + {1\over 4}\tau.\]
Here $\tau$ denotes the scalar curvature of the Riemannian
manifold $M$.
\end{theorem}

\Proof: We only have to prove that ${\cal R} = {1\over 4}\tau$. It suffices 
to show that with respect to a local orthonormal frame $\{E_1,\dots,E_m\}$ of 
$TM$ the curvature $\curv (\nabla^S)$ is given by 
\[\curv(\nabla^S)(X,Y) = {1\over 4}\sum_{k,\ell = 1}^m \l
R(X,Y)E_k,E_\ell\r E_k \cdot E_\ell,\quad X,Y \in
C^\infty(TM\mid_U)\hfill(\ast)\]
since then we obtain
\begin{eqnarray*}\hspace{-1cm}
{\cal R}\mkern-30mu&&={1\over 2}\sum_{i,j = 1}^m E_i \cdot E_j \cdot
\curv(\nabla^S)(E_i,E_j)\\&&={1\over 8}\sum_{i,j,k,\ell = 1}^m \l
R(E_i,E_j)E_k,E_\ell\r E_i \cdot E_j \cdot E_k \cdot E_\ell\\ &&={1\over
8}\sum_{\ell = 1}^m\Big(\mkern-5mu\sum_{i\ne j \ne k \ne \ell}\mkern-10mu\l
R(E_i,E_j)E_k + R(E_k,E_i)E_j + R(E_j,E_k)E_i,E_\ell\r E_i \cdot E_j \cdot
E_k\\ &&\quad +\sum_{i,j}\l R(E_i,E_j)E_i,E_\ell\r E_i\cdot E_j\cdot E_i
+\sum_{i,j}\l R(E_i,E_j)E_j,E_\ell\r E_i\cdot E_j\cdot E_j\Big)E_\ell\\
&&=-{1\over 4}\sum_{i,j,\ell = 1}^m \l R(E_i,E_j)E_j,E_\ell\r E_i\cdot E_\ell
 = {1\over 4}\sum_{i,j = 1}^m \l R(E_i,E_j)E_j,E_i\r= {1\over 4}\tau
\end{eqnarray*}
by the symmetries of the Riemann curvature tensor $R$ and by the
definition of $\tau$.

Now it is a straight-forward computation to show that for fixed
vector fields $X$ and $Y$ the right-hand side of $(\ast)$ which we
denote by $R(X,Y)$ shares the same properties (R) and (D) as the
left hand-side and so does their difference
$T=\curv(\nabla)(X,Y)-R(X,Y)$. In particular, by (D) it commutes
with the left-action of $\Cl M$ and so acts as multiplication by
an element $\gamma\in C^\infty(M,\C)$ which by (R) is
skew-adjoint, i.e.\ $\gamma=i\eta$ with $\eta\in C^\infty(M,\R)$.
Actually, $\eta$ has to vanish, since $T$ also respects the real
structure on $S$, i.e.\ commutes with the structural map $\theta$.
\endproof
\medskip

This Bochner-Weitzenb\"{o}ck formula for the Spin-Dirac operator is used by A.\ 
Lichnerowicz [Lic] to prove the following vanishing theorem. The relation of 
$D_{AS}^2$ and the scalar curvature had however already been noted by E.\ 
Schr\"odinger in 1932 [Sch]. 
\medskip

\noindent 
{\bf Corollary}~ {\rm (Lichnerowicz)}\quad {\it Let $M$ be a 
compact spin manifold with positive scalar curvature. Then there are no 
harmonic spinors on $M$. If $\dim M=4k$, then $\hat A(M)=0$.} 
\smallskip

\Proof: The first assertion is immediate while the second one is a
consequence of the Atiyah-Singer index theorem.
\endproof
\medskip

We also study the twisted Dirac operator $\nd \otimes I_E$, where $E$ is a 
Hermitian vector bundle $E$ with connection $\nabla^E$, i.e.\ the Dirac 
operator of the Dirac bundle $(S \otimes E,\nabla^{S \otimes E})$. Let ${\cal 
R}^E : C^\infty(S \otimes E) \to C^\infty(S \otimes E)$ denote the zero order 
differential operator, which for sections $\sigma\otimes \fs$ and the frame 
$(E_i)_{1 \le i \le m}$ is defined by 
\[{\cal R}^E(\sigma \otimes
\fs) = {1\over 2}\sum_{j,k=1}^m E_j\cdot E_k\cdot\sigma
\otimes\curv(\nabla^E)(E_j,E_k)\fs.\]

\begin{theorem} Let $M$ be spin and $S$ and $E$ as
before. Then the Spin-Dirac operator $D_{AS} \otimes I_E$ and the
Bochner-Laplace operator $\nabla^\ast\nabla$ of the tensor bundle
$S\otimes E$ are related by
\[(D_{AS} \otimes I_E)^2 = \nabla^\ast
\nabla + {1\over 4}\tau + {\cal R}^E.\]
Here $\tau$ is again the
scalar curvature of $M$ .
\end{theorem}

\Proof: For $\sigma \in C^\infty(S)$ and $\fs \in C^\infty(E)$ we
have
\[\nabla^{S \otimes E}(\sigma \otimes \fs) = (\nabla^S
\sigma)\otimes \fs + \sigma (\nabla^E\fs).\]
This entails
\[\curv(\nabla^{S\otimes E})(\sigma\otimes \fs) =
\curv(\nabla^S)(\sigma)\otimes \fs + \sigma \otimes
\curv(\nabla^E)(\fs)\]
and
\begin{eqnarray*}
{\cal R}(\sigma \otimes \fs) &=& {1\over 2}\sum_{j,k=1}^mE_j \cdot
E_k\cdot\curv(\nabla^{S\otimes E})(E_j,E_k)(\sigma \otimes \fs)\\
&=& {1\over 2}\sum_{j,k=1}^m E_j \cdot
E_k\cdot\big(\curv(\nabla^{S})(E_j,E_k)(\sigma)\big) \otimes \fs\\
&\quad\quad&+{1\over 2}\sum_{j,k=1}^m E_j \cdot
E_k\cdot\sigma\otimes\curv(\nabla^{ E})(E_j,E_k)(\fs)\\ &=&
{1\over 4}\tau(\sigma\otimes \fs) + {\cal R}^E(\sigma \otimes
\fs).
\end{eqnarray*}

\vspace{-12pt}\rightline\qed
\medskip

{\bf Remark}\quad For the Spin$^c$-Dirac operator $D_{S^c}$ there
is a Bochner-Weit\-zen\-b\"{o}ck formula, too. If $M$ is spin a
spin$^c$ structure is given by $S^c=S\otimes L$ for some complex
line bundle $L$. Choosing the product connection on $S^c$ with
some Hermitian connection $\nabla^L$ on $L$ the square of the
corresponding Spin$^c$-Dirac operator $D_{S^c}$ satisfies
\[D_{S^c}^2 = \nabla^\ast \nabla + {1\over 4}\tau + {\cal R}^L.\]
Because of
\begin{eqnarray*}
{\cal R}^L(\sigma\otimes\fs)&=&{1\over 2}\sum_{j,k=1}^mE_j \cdot
E_k\cdot\sigma\otimes\curv(\nabla^{L})(E_j,E_k)(\fs)\\ &=&{1\over
2}\sum_{j,k=1}^mE_j \cdot E_k\cdot\sigma\otimes\Omega^L(E_j,E_k)(\fs)\\
&=&{1\over 2}\sum_{j,k=1}^m\Omega^L(E_j,E_k)E_j \cdot
E_k\cdot\sigma\otimes\fs=(\Omega^L\cdot\sigma)\otimes\fs
\end{eqnarray*}
(with $\Omega^L$ denoting the curvature form of $\nabla^L$) we
obtain
\[D_{S^c}^2 = \nabla^\ast \nabla + {1\over 4}\tau + \Omega^L.\]
Since this computation is local, we can also apply it in the
general non-spin case. Although $S^c$ is a product $S\otimes L$
only locally the line bundle $L_{S^c}=\Hom_{\Cl
M_\C}(\o{S^c},S^c)= L\otimes L$ is nevertheless globally defined.
Choosing a Hermitian connection, the corresponding curvature form
$\Omega$ satisfies $\Omega=2\Omega^L$. Thus we obtain the
Weitzenb\"{o}ck formula
\[D_{S^c}^2 = \nabla^\ast \nabla + {1\over 4}\tau + {1\over 2}\Omega.\]
The index formula for $D^0_{S^c}$ can also be established by a
local computation (cf.\ [Sdr2]). With $c=c_1(L_{S^c})$ one obtains
\[\ind D^0_{S^c}=\int_M e^{c/2}\hat A(TM).\]

The non-vanishing of the $\hat A$-genus is the simplest
obstruction for a Riemannian metric with positive scalar
curvature. N.\ Hitchin [Hit] has introduced an invariant
$\alpha(M)$, which can be defined for spin manifolds $M$ of any
dimension and which coincides with $\hat A(M)$ if $m=4k$. It again
vanishes in case of positive scalar curvature. For simply
connected manifolds $\alpha(M)=0$ is even sufficient for such a
metric to exist as S.\ Stolz [Sto] proved in 1989; cf.\ [RS] for a
survey of the current state.\\The Bochner-Weitzenb\"{o}ck formula for
the Spin$^c$-Dirac operator on oriented compact 4-manifolds is the
footing of the so-called Seiberg-Witten theory in which the
theoretical physicists N.\ Seiberg and E.\ Witten initiated new
differential topological invariants in 1994. These lead to new
essential contributions for the classification of 4-manifolds
[Mor].
\bigskip
\bigskip

{\bf\Large References} 

\bigskip

[AT] Anderson, F.W., Fuller, K.R.: {\it Rings and Categories of
Modules,} Springer, New York - Heidelberg - Berlin, 1974\smallskip

[ABS] Atiyah, M.F., Bott, R., Shapiro, A.: Clifford modules, {\it
Topology} 3 (1964), Suppl.\ 1, 3-38\smallskip

[AS] Atiyah, M.F., Singer, I.M.: The index of elliptic operators
on compact manifolds, {\it Bull.\ Amer.\ Math.\ Soc.} 69 (1963)
422-432\smallskip

[BD] Baum, P., Douglas, R.G.: Index theory, bordism, and
$K$-homology, in: Operator Algebras and $K$-Theory, {\it Contemp.\
Math.} 10, 1-31, Amer.\ Math.\ Soc., Providence, RI,
1982\smallskip

[Boc] Bochner, S.: Vector fields and Ricci curvature, {\it Bull.\
Amer.\ Math.\ Soc.} 52 (1946) 776-797\smallskip

[BBW] Booss-Bavnbek, B., Wojciechowski, K.P.: {\it Elliptic
Boundary Problems for Dirac Operators}, Birkh\"{a}user, Basel-Boston,
1993\smallskip

[BH] Borel A., Hirzebruch, F.: Characteristic classes on
homogeneous soaces II, {\it Amer.\ J.\ Math.} 81 (1959)
315-382\smallskip

[BW] Brauer, R., Weyl, H.: Spinors in $n$ dimensions, {\it Amer.\
J.\ Math.} 57 (1935) 425-449\smallskip

[Car] Cartan, E.: Sur les groupes projectifs qui laissent
invariante aucune multiplicit\'{e} plane, {\it Bull.\ Soc.\ Math.\
France} 41 (1913) 53-96\smallskip

[Che] Chevalley, C.: {\it The Algebraic Theory of Spinors},
Columbia Univ.\ Press, New York, 1954, together with {\it The
Construction and Study of Certain Important Algebras}, Math.\
Soc.\ Japan, Tokyo, 1955, reprinted in: Collected Works Vol.\ 2,
Springer, New York - Heidelberg - Berlin, 1998\smallskip

[Cli] Clifford, W.K.: Applications of Grassmann's extensive
algebra, {\it Amer.\ J.\ Math.} 1 (1878) 350-358 or {\it Math.\
Papers}, 266-276; cf.\ also: On the classification of geometric
algebras [1876], ibid.\ 397-401\smallskip

[Con1] Connes, A.: {\it Noncommutative Geometry}, Academic Press,
New York, 1994\smallskip

[Con2] Connes, A.: Noncommutative geometry and reality, {\it J.\
Math.\ Phys.} 36 (11) (1995) 6194-6231\smallskip

[Con3] Connes, A.: Gravity coupled with matter and foundation of
noncommutative geometry, {\it Commun.\ Math.\ Phys.} 182 (1996)
155-176\smallskip

[Dar] Darwin, C.G.: The wave equation of the electron, {\it Proc.\
Roy.\ Soc.\ London} (A) 118 (1928) 654-680\smallskip

[Dir] Dirac, P.A.M.: The quantum theory of the electron, {\it
Proc.\ Roy.\ Soc.\ London} (A) 117 (1927) 610-624\smallskip

[Eck] Eckmann, B.: Hurwitz-Radon matrices revisited: from
effective solution of the Hurwitz matrix equations to Bott
periodicity, {\it CRM Proc.\ \& Lecture Notes} 6, 23-35, Amer.\
Math.\ Soc., 1994\smallskip

[Fri] Friedrich, T.: {\it Dirac Operatoren in der Riemannschen Geometrie},
Vieweg, Braunschweig, 1997\smallskip

[Gil] Gilkey, P.B.: {\it Invariance Theory, the Heat Equation, and
the Atiyah-Singer Index Theorem,} (2nd ed.), CRC, Baton Rouge,
1995\smallskip

[GH] Greub, W.H., Halperin, S.: An intrinsic definition of the
Dirac operator, {\it Collect.\ Math.} 26 (1975) 19-37\smallskip

[Hae] Haefliger, A.: Sur l'extension du groupe structural d'un
espace fibr\'{e}, {\it C.\ R.\ Acad.\ Sci.\ Paris} 243 (1956)
558-560\smallskip

[Hit] Hitchin, N.: Harmonic spinors, {\it Adv.\ Math.} 14 (1974)
1-55\smallskip

[HH] Hirzebruch, F., Hopf, H.: Felder von Fl\"{a}chenelementen in
4-dimensionalen Mannigfaltigkeiten, {\it Math.\ Ann.} 136 (1958)
156-172\smallskip

[Hur] Hurwitz, A.: \"{U}ber die Komposition der quadratischen Formen,
{Math.\ Ann.} 88 (1923) 1-25 (posthumous)\smallskip

[JW] Jordan, P., Wigner, E.: \"{U}ber das Paulische \"{A}quivalenzverbot,
{\it Zeit.\ f.\ Physik} 47 (1927) 631-651\smallskip

[Kae] Kaehler, E.: Der innere Differentialkalk\"{u}l, {\it Rend.\
Mat.} 21 (1962) 425-523\smallskip

[Krb] Karoubi, M.: {\it $K$-Theory,} Springer, New York -
Heidelberg - Berlin, 1974\smallskip

[Kar1] Karrer, G.: Einf\"{u}hrung von Spinoren auf Riemannschen
Mannigfaltigkeiten, {\it Ann.\ Acad.\ Scient.\ Fennicae} Ser.\ A.\
I.\ Mathematica 336/5 (1963)\smallskip

[Kar2] Karrer, G.: Darstellung von Cliffordb\"{u}ndeln, {\it Ann.\
Acad.\ Scient.\ Fennicae} Ser.\ A.\ I.\ Mathematica 521
(1973)\smallskip

[LM] Lawson, H.B., Michelsohn, M.L.: {\it Spin Geometry,}
Princeton Univ.\ Press, Princeton, 1989\smallskip

[Lic] Lichnerowicz, A.: Spineurs harmoniques, {\it C.\ R.\ Acad.\
Sci.\ Paris} A 257 (1963) 7-9\smallskip

[Lip] Lipschitz, R.: {\it Untersuchungen ueber die Summen von
Quadraten}, Max Cohen \& S., Bonn, 1886 (French extract in {\it
Bull.\ Sci.\ Math.\ 2 S\'{e}r.} 10 (1886) 163-183)\smallskip

[Mil] Milnor, J.: Spin structures on manifolds, {\it l'Ens.\
Math.} 9 (1963) 198-203\smallskip

[Mor] Morgan, J.: {\it The Seiberg-Witten Equations and
Applictions to the Topology of Smooth 4-Manifolds}, Princeton
Univ.\ Press, Princeton, 1996\smallskip

[Pau] Pauli, W.: Zur Quantenmechanik des magnetischen Elektrons,
{\it Zeit.\ f.\ Physik} 43 (1927) 601-623\smallskip

[Ply] Plymen, R.: Strong Morita equivalence, spinors and
symplectic spinors, {\it J.\ Operator Theory} 16 (1986)
305-324\smallskip

[Rad] Radon J.: Lineare Scharen orthogonaler Matrizen, {\it Abh.\
Math.\ Seminar Hamburg} 1 (1922) 1-14\smallskip

[Ren] Rennie, A.: Commutative geometries are spin manifolds,
preprint, Univ. Adelaide,
http://xxx.lanl.gov/math-ph/9903021\smallskip

[RS] Rosenberg, J.M., Stolz, S.: Manifolds of positive scalar
curvature, {\it Algebraic Topology and its Applications}, 241-267,
MSRI Publ.\ 27, Springer, New York, 1994\smallskip

[Sdr1] Schr\"{o}der, H.: Funktionalanalysis, Verlag Harri Deutsch,
Frankfurt a.M., 2000\smallskip

[Sdr2] Schr\"{o}der, H.: Globale Analysis, Textbook (manuscript),
Univ.\ Dortmund, 2000\smallskip

[Sch] Schr\"{o}dinger, E.: Diracsches Elektron im Schwerefeld I., {\it
Sitz.-ber.\ Preuss.\ Akad.\ Wiss.\ Berlin, Phys.-Math.\ Kl.} XI
(1932) 105-128\smallskip

[Sto] Stolz, S.: Simply connected manifolds of positive scalar
curvature, {\it Bull.\ Amer.\ Math.\ Soc.} 23 (1990)
427-432\smallskip

[Var] V\'{a}rilly, J.C.: {\it An Introduction to Noncommutative
Geometry}, Lecture Notes EMS Summer School on Noncommutative
Geom.\ and Appl., Monsaraz and Lisabon, 1997,
http://xxx.lanl.gov/physics/9709045\smallskip

[vdW1] Waerden, B.\ L.\ van der: Exclusion principle and spin, in:
{\it Theoretical Physics in the Twentieth Century}, (eds. M.\
Fierz, V.F.\ Weisskopf), Interscience Publ., New York, 1960, pp.
199-244\smallskip

[vdW2] Waerden, B.\ L.\ van der: On Clifford algebras, {\it
Indag.\ Math.} 28 (1966) 78-83\smallskip

[Wei] Weil, A.: Correspondence, {\it Ann.\ of Math.} 69 (1959)
247-251\smallskip

[Wey1] Weyl, H.: {\it The Theory of Groups and Quantum Mechanics},
Dover Publ., New York, (transl.\ by H.P.\ Robertson from {\it
Gruppentheorie und Quantenmechanik}, Hirzel, Leipzig,
1931)\smallskip

[Wey2] Weyl, H.: {\it The classical groups}, Princeton Univ.\
Press, Princeton, 1949$^2$\smallskip

\vspace{1cm}

Address:\\ 
Herbert Schr\"{o}der\\ 
Fachbereich Mathematik\\ 
Universit\"{a}t Dortmund\\ 
Postfach 50 05 00\\ 
D-44221 Dortmund \\ 
e-mail: schroed@math.uni-dortmund.de 
 
\end{document}